\RequirePackage{expl3}
\documentclass{ieeeaccess}
\usepackage{cite}
\usepackage{amsmath,amssymb,amsfonts}
\usepackage{algorithmic}
\usepackage{graphicx}
\usepackage{textcomp}

\usepackage{mathptmx} 
\usepackage{times} 
\usepackage{siunitx}
\usepackage{tabularx,booktabs}
\usepackage{bm}
\usepackage[whole]{bxcjkjatype}
\usepackage{color}
\usepackage[table]{xcolor}
\usepackage{subcaption}
\usepackage{mathtools}
\usepackage{dblfloatfix} 

\newcommand{\norm}[1]{\lVert #1 \rVert}

\newcommand{\abs}[1]{\lvert #1 \rvert}

\usepackage{accents}

\DeclareMathOperator{\diag}{diag}

\DeclareMathOperator{\asym}{asym}
\DeclareMathOperator{\stat}{stat}

\DeclareSymbolFont{bbold}{U}{bbold}{m}{n}
\DeclareSymbolFontAlphabet{\mathbbold}{bbold}
\newcommand{\onev}{\mathbbold{1}}

\usepackage{calrsfs}
\DeclareMathAlphabet{\pazocal}{OMS}{zplm}{m}{n}
\renewcommand{\mathcal}[1]{\pazocal{#1}}

\newtheorem{definition}{Definition}[section]

\newtheorem{lemma}[definition]{Lemma}
\newtheorem{proposition}[definition]{Proposition}

\newtheorem{remark}[definition]{Remark}

\newtheorem{problem}[definition]{Problem}

\def\BibTeX{{\rm B\kern-.05em{\sc i\kern-.025em b}\kern-.08em
    T\kern-.1667em\lower.7ex\hbox{E}\kern-.125emX}}
\begin{document}
\history{Date of publication xxxx 00, 0000, date of current version xxxx 00, 0000.}
\doi{10.1109/ACCESS.2017.DOI}

\graphicspath{{../}}

\title{Deep Learning-Based Average Consensus}
\author{\uppercase{Masako Kishida}\authorrefmark{1}, \IEEEmembership{Senior Member, IEEE},
  \uppercase{Masaki Ogura\authorrefmark{2}, \IEEEmembership{Member, IEEE},
    Yuichi Yoshida\authorrefmark{1} and Tadashi Wadayama\authorrefmark{3}, \IEEEmembership{Member, IEEE}}}
\address[1]{Principles of Informatics Research Division, National Institute of Informatics, Tokyo
  101-8430, Japan (e-mail: kishida, yyoshida@nii.ac.jp)}
\address[2]{Graduate School of Information Science and Technology, Osaka University, Osaka 565-0871, Japan (e-mail: m-ogura@ist.osaka-u.ac.jp)}
\address[3]{Nagoya Institute of Technology, Nagoya, Aichi 465-8555, Japan (e-mail: wadayama@nitech.ac.jp)}
\tfootnote{M.K. and M.O. contributed equally to this work.\\
M.O.~is supported by The Telecommunications Advancement Foundation. Y.Y. is supported by JSPS KAKENHI Grant Number 18H05291.
}
%\tfootnote{This paragraph of the first footnote will contain support 
%information, including sponsor and financial support acknowledgment. For 
%example, ``This work was supported in part by the U.S. Department of 
%Commerce under Grant BS123456.''}

\markboth
{Author \headeretal: Preparation of Papers for IEEE TRANSACTIONS and JOURNALS}
{Author \headeretal: Preparation of Papers for IEEE TRANSACTIONS and JOURNALS}

\corresp{Corresponding authors: Masako Kishida (e-mail: kishida@nii.ac.jp) and Masaki Ogura (e-mail: m-ogura@ist.osaka-u.ac.jp).}

\begin{abstract}
In this study, we analyzed the problem of accelerating the linear average consensus algorithm for complex networks. We propose a data-driven approach to tuning the weights of temporal (i.e., time-varying) networks using deep learning techniques. Given a finite-time window, the proposed approach first unfolds the linear average consensus protocol to obtain a feedforward signal-flow graph, which is regarded as a neural network. The edge weights of the obtained neural network are then trained using standard deep learning techniques to minimize consensus error over a given finite-time window. Through this training process, we obtain a set of optimized time-varying weights, which yield faster consensus for a complex network. We also demonstrate that the proposed approach can be extended for infinite-time window problems. Numerical experiments revealed that our approach can achieve a significantly smaller consensus error compared to baseline strategies.
\end{abstract}


\titlepgskip=-15pt

\maketitle

\section{Introduction}

The distributed agreement problem for networks, which is often referred to as the consensus problem~\cite{Olfati-Saber2007}, is an important problem in network science and engineering with applications in load balancing~\cite{Cyb89}, data fusion~\cite{XiaBL05}, multi-agent coordination~\cite{Ren2005}, distributed computing~\cite{Xiao2004}, distributed sensor networks~\cite{Cortes2005}, wireless communication systems~\cite{Senel2017}, and power systems~\cite{Dorfler2012}. Recently, this problem has also appeared in online machine learning procedures for processing big data (see \cite{CheS12,TsiLR12} and the references therein).

In the \emph{average} consensus problem, the nodes in a network seek to converge their states to the average of their initial states in a distributed manner. The standard solution to this problem is to use the linear average consensus algorithm~\cite{Olfati-Saber2004}, where each node updates its state by calculating the weighted linear average of its own state and the states of its neighbors. This algorithm generates a linear dynamical system whose state transition matrix involves the Laplacian matrix of the underlying communication network.

Designing consensus algorithms with fast convergence is of significant practical importance because such algorithms allow multi-agent systems to reach agreements in fewer iterations, meaning they will consume less communication resource. In the context of the linear average consensus algorithm, the problem of finding the optimal weights of edges to maximize asymptotic consensus speed can be reduced to a convex optimization problem~\cite{Xiao2004} under the assumption that the communication network is static and undirected. The authors of \cite{Kempton2018} recently demonstrated that optimal weights can be computed in a distributed manner via iterative computations. Zelazo et al.~\cite{Zelazo2013} clarified the role of cycles in the linear average consensus algorithm and proposed an approach for accelerating this algorithm by adding new edges to a network. For the case of directed networks, Hao and Barooah~\cite{Hao2012} proposed a method to accelerate the convergence rate of a linear (but not necessarily average) consensus algorithm by tuning the weights of edges in the target network. Additionally, for the case of gossip algorithms, the authors of \cite{BoyGP06} demonstrated that it is possible to tune weights in a distributed manner such that the weights converge to the optimal weights, which yields solutions for averaging problems.

A natural consequence of seeking further acceleration of consensus algorithms is the emergence of finite-time consensus algorithms~\cite{Sundaram2007}, where edge weights are typically assumed to vary with time and designers exploit this additional flexibility to realize consensus in finite time. The finite-time consensus algorithm proposed in \cite{SanKM14} achieves a consensus using a time-invariant weight updating law. The algorithm proposed in \cite{Hendrickx2015} achieves consensus using stochastic (possibly asymmetric) matrices in $N(N-1)/2$ iterations, where $N$ denotes the number of nodes in the target network. The authors of \cite{Safavi2015,Shang2016} used graph signal processing tools (e.g., \cite{Shuman2013}) to demonstrate that by allowing non-stochasticity for state-updating matrices, one can realize finite-time consensus in at most $N$ steps. The theoretical aspects of these works were further investigated in~\cite{Apers2017}. Recently, the authors of \cite{Falsone2018} showed that the number of steps required for consensus can be further improved to $N/2$ in the specific case of a ring network with an even number of nodes.

Despite the advances in consensus acceleration described above, there is still a need for an effective approach to answering the following basic question. Given a finite-time window and underlying network structure, how should one dynamically tune the edge weights in the network to achieve the most accurate consensus possible within a specific time window? If the length of a time window is insufficient for executing the aforementioned finite-time consensus algorithms, then currently available options are effectively limited to using static optimal strategies (e.g., \cite{Xiao2004}), which do not allow one to tune the weights of a network dynamically. To fill this gap, we propose a data-driven approach to tuning the weights of undirected temporal (i.e., time varying) networks using deep learning techniques. We first unfold the consensus algorithm to obtain a feedforward signal-flow graph \cite{Ito2019}, which we regard as a neural network. We then use the standard stochastic gradient descent algorithm to update the parameters in each layer of this neural network (i.e., the weights of each snapshot of the temporal network) to minimize consensus error over a finite-time window, which yields an optimized temporal network for faster consensus. We numerically confirm that our approach can significantly accelerate the convergence speed of the linear average consensus algorithm.

The remainder of this paper is organized as follows. In Section~\ref{sec:problemSetting}, we define the problem of dynamically tuning edge weights to accelerate the linear average consensus algorithm and present our approach to solving this problem using standard techniques from the field of deep learning. In Section~\ref{sec:numericalSimulations}, we evaluate the performance of the proposed method through various numerical experiments. We conclude this paper in Section~\ref{sec:conc}.

\section{Weight Optimization using Deep Learning Techniques} \label{sec:problemSetting}

In this section, we describe our approach to tuning the edge weights of a network to accelerate the linear average consensus algorithm within a given finite-time window. We first provide a brief review of the linear average consensus algorithm and state its basic properties. We then describe our data-driven approach to tuning the weights of a network, in which we apply techniques from the deep learning to the signal-flow graph obtained by unfolding the consensus algorithm.

\subsection{Linear average consensus algorithm}

Let $G$ be an undirected, unweighted, and connected network with a node set~$V= \{1, \dotsc, N\}$ and edge set $E$, which consists of unordered pairs of nodes in~$V$. Each node in $G$ represents an agent that is supposed to communicate with its neighbors at each time step. In this paper, we focus on discrete-time dynamics. Let $x_i(k) \in \mathbb{R}$ denote the state of the $i$th node at time~$k\geq 0$ and let $\mathcal N_i$ denote the set of neighbors of node~$i$. In the standard linear average consensus protocol~\cite{Olfati-Saber2007}, each node~$i$ updates its own state according to the following difference equation:
\begin{equation}\label{eq:weightedTimevaryingConsensus}
  \!\!x_i(k+1) \!=\! x_i(k)\! +\! \sum_{j\in \mathcal N_i} w_{ij}(k)(x_j(k) - x_i(k)), \
  x_i(0) =  x_{0,i},
\end{equation}
where $w_{ij}(k) = w_{ji}(k)\geq 0$ represents the weight of the (undirected) edge~$\{i, j\}$ at time~$k$
and $x_{0,i}$ denotes the initial state of node~$i$. Although we assume the symmetry of edges for convenience, our framework can be easily extended to asymmetric cases.
At each time~$k\geq 0$, we define the $(i,j)$ element of the adjacency matrix of the network $W(k)\in \mathbb{R}^{N\times N}$ as
\begin{equation}
  W_{ij}(k) = \begin{cases}
    w_{ij}(k), & \mbox{if $j \in \mathcal N_i$},
    \\
    0,         & \mbox{otherwise}
  \end{cases} \label{eq:W}
\end{equation}
and the degree matrix of the network at time~$k$ as
\begin{equation}
  D(k) = \diag(d_1(k), \dots, d_N(k)),\ \ d_i(k) = \sum_{j\in\mathcal N_i} w_{ij}(k).
\end{equation}
Then, by using the Laplacian matrix of the network
\begin{equation}
  L(k) = D(k)-W(k),
\end{equation}
the evolution of the state vector
\begin{equation}
  x(k) = \begin{bmatrix}
    x_1(k) & \cdots & x_N(k)
  \end{bmatrix}^\top
\end{equation}
in the linear average consensus protocol (\ref{eq:weightedTimevaryingConsensus}) can be written as
\begin{equation}
  x(k+1) = (I-L(k))x(k),
  \ \
  x(0) = x_0, \label{eq:protocol}
\end{equation}
where
$
  x_0 = \begin{bmatrix}
    x_{1, 0} & \cdots & x_{N, 0}
  \end{bmatrix}^\top $
denotes the initial state vector.

\begin{figure}[tb] 
  \centering
    \includegraphics[width=1\linewidth]{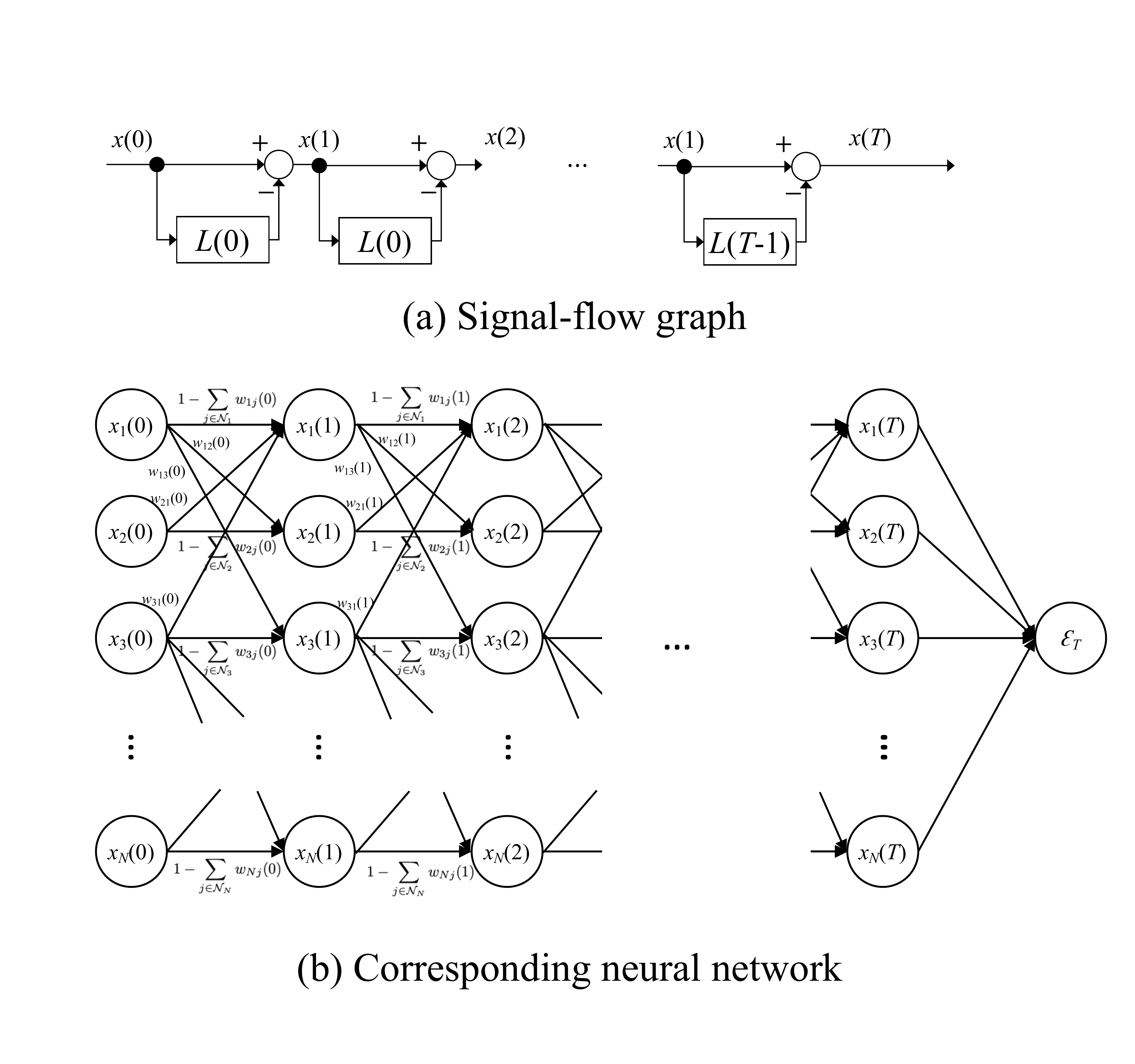}
    \caption{Unfolded signal-flow graph and corresponding neural network. Signal-flow graph was obtained by unfolding the linear average consensus algorithm~(\ref{eq:weightedTimevaryingConsensus}) for $k=0, \dotsc, T-1$. \label{fig:NeuralNet}}
\end{figure}

%\begin{figure*}[tb]
%  \centering
%  \begin{subfigure}{.35\linewidth}
%    \centering
%    \includegraphics[width=1\linewidth]{NN.pdf}
%    \caption{Signal-flow graph obtained by unfolding the linear average consensus algorithm~(\ref{eq:weightedTimevaryingConsensus}) for $k=0, \dotsc, T-1$.}
%    \label{fig:NeuralNet}
%  \end{subfigure}
%  \hfil
%  \begin{subfigure}{.575\linewidth}
%    \includegraphics[width=.8\linewidth]{NN2.pdf}
%    \caption{Corresponding neural network}
%    \label{fig:NeuralNet2}
%  \end{subfigure}
%  \caption{Unfolded signal-flow graph and corresponding neural network}
%\end{figure*}

The objective of this paper is to present a framework for tuning the weights $\{w_{ij}(k)\}_{k\geq 0, \{i, j\} \in E}$ to achieve  a faster average consensus within a given finite-time window. We denote the average of the initial states of the nodes as
\begin{equation*}
  c = \frac{1}{N}\sum_{i=1}^N x_{i,0}.
\end{equation*}
Define the consensus error vector
\begin{equation*}
  e(k) = x(k) - c \onev,
\end{equation*}
where $\onev$ denotes an all-ones $N$-dimensional column vector. We are now ready to formally state the problem studied in this paper.

\begin{problem}[Consensus acceleration problem]\label{prb:}
Let $G$ be an undirected and unweighted network containing $N$ nodes. Let $T$ be a given positive integer. Assume that the set of initial states follows a probability distribution $\mathcal X_0$, i.e.,
\begin{equation*}
  \{x_{0, 1}, \dotsc, x_{0, N} \}\sim \mathcal X_0.
\end{equation*}
Find the set of nonnegative weights
\begin{equation*}
  \{w_{ij}(k)\}_{k \in \{0, 1, \cdots,  T-1\}, \{i, j\} \in E}
\end{equation*}
that minimizes the average consensus error defined by
\begin{equation*}
  \varepsilon_{\,T} = \mathbb{E} [\norm{e(T)} ],
\end{equation*}
where $\lVert \cdot \rVert$ denotes the Euclidean norm in $\mathbb{R}^N$, and $\mathbb{E}[\cdot]$ denotes the expected value.
\end{problem}

Because Problem \ref{prb:} is a non-convex problem, it is difficult to compute a set $\{w_{ij}(k)\}_{k\geq 0, \{i, j\} \in E}$ that minimizes $\varepsilon_{\,T}$. This difficulty motivated us to tackle this problem by using a data-driven approach to find a suboptimal solution.
In the next subsection, we describe our data-driven approach to tuning edge weights using deep learning techniques. It should be noted that although we assume some knowledge of the initial probability distribution~$\mathcal X_0$ during the process of optimization, optimized edge weights can significantly accelerate the consensus protocol even if the initial states do not follow the assumed distribution. We numerically demonstrate this universal property of our approach in Subsection~\ref{subsec:deterministic}.

\subsection{Data-driven weight optimization}\label{subsec:optimization}

To adjust edge weights using deep learning techniques, we first unfold the recursive state-update formula~(\ref{eq:weightedTimevaryingConsensus}) to obtain a signal-flow graph, as shown in Fig.~\ref{fig:NeuralNet}(a). Unlike a standard deep neural network, the resulting neural network has a fixed structure and uses no activation functions. The structure between two layers of this neural network corresponds to the structure of the graph $G$ and is the same between all layers. The neurons in the $k$th layer correspond to the nodes at time $k$, as shown in Fig.~\ref{fig:NeuralNet}(b).

We then apply a standard technique from the field of deep learning to adjust the weights.
We use the mean squared error, $\varepsilon_{\,T}^2$, as a loss function. We do not append any regularization terms.
As in~\cite{Ito2019}, we use the technique of the incremental training to adjust the weights during the training process. For incremental training, we first consider only the first layer (i.e., we set $k=1$ in Fig.~\ref{fig:NeuralNet}(a)) and attempt to minimize the loss function of the average consensus error $\varepsilon_{\,1}^2$ using a number of randomly generated initial states $x_0$ as the training data. This first step is called the first generation.  After training the first set of weights $w_{ij}(0)$, we proceed to train the first two sets of edge weights by appending the second layer to the neural network and replacing the loss function with $\varepsilon_{\,2}^2$. For this training, we use the result from the 1st generation as the initial values for the first layer and train the entire neural network. We repeat this process to optimize the all of the weights $w_{ij}(T-1)$ between ${T-1}$st and $T$th layers by minimizing  $\varepsilon_{\,T}^2$.

Therefore, using the set of equations \eqref{eq:W}-\eqref{eq:protocol} with the optimized weights, we can obtain the optimized weighted Laplacian matrices for the network
$L^\star(0)$, \dots, $L^\star(T-1)$ for given $G$ and $T$. The proposed consensus algorithm is defined as
\begin{equation}
  x(k+1) = (I-L^\star(k))x(k),\quad {0\leq k\leq T-1}. \label{eq:protocol2}
\end{equation}

The proposed approach for solving Problem \ref{prb:} can be considered as supervised learning with an input object $x_0$ and desired output value $c \onev$, where $c = \onev^\top x_0/N$, with identity activation functions over a neural network with a specified structure.

\subsection{Periodic continuation} \label{sec:periodic}

For very long or infinite-time windows, where it is not necessarily feasible to tune weights using the proposed approach, we may adopt a periodic continuation of the proposed algorithm.
By periodically extending the state transition matrices $I-L^\star(0)$, \dots, $I-L^\star(T-1)$ obtained previously, the average consensus protocol for a time greater than $T$ can be defined as follows:
\begin{align}\begin{aligned}\label{eq:periodicConsensnsProt}
    x(sT+\tau+1) & =\left(\prod_{t=0}^{\tau}(I-L^\star(\tau-t)) \right)x(sT), \\
                 & \ 0\leq \tau \leq T-1, \ s\geq 0.
  \end{aligned}
\end{align}

To assess the performance of this algorithm, we introduce the following quantity:
\begin{equation}\label{eq:asymConvFact}
  r_{\asym} =
  \adjustlimits\sup_{x_0\neq c\onev} \limsup_{k\to\infty}
  \left( \frac{\norm{e(k)}}{\norm{e(0)}}\right)^{1/k}, 
\end{equation}
which we call the asymptotic convergence factor. The next lemma provides an explicit representation of the asymptotic convergence factor \eqref{eq:asymConvFact} for the consensus algorithm~(\ref{eq:periodicConsensnsProt}).

\begin{proposition}\label{LEM:}
Given $G$ and $T$, let $L^\star(0)$, \dots, $L^\star(T-1)$ denote the optimized weighted Laplacian matrices of the networks derived by our deep learning algorithm.
  Define a matrix
    \begin{equation*}
      M = \prod_{t=0}^{T-1}(I-L^\star(T-1-t))
    \end{equation*}
    and let $\sigma(M)$ denote the set of eigenvalues of $M$. Additionally, define
    \begin{equation}\label{eq:asmpConvFactorDL}
      \rho = \max\{ \abs{\lambda}\mid \lambda \in \sigma(M) \backslash \{1\}\}.
    \end{equation}
   Assume that the eigenvalue $1$ of $M$ is simple. Then, the asymptotic convergence factor of the periodic continuation with a period $T$ given by the consensus algorithm~(\ref{eq:periodicConsensnsProt}) is equal to $\rho^{1/T}$.
\end{proposition}

\begin{remark}
It should be noted that
  \begin{equation*}
       r_{\asym}^T =\! \! \adjustlimits\sup_{x_0\neq c\onev} \limsup_{k\to\infty} \left(\dfrac{\norm{e(k)}}{\norm{e(0)}}\right)^{1/k}
    \end{equation*}
   is  the asymptotic convergence factor of the periodic continuation with a period $T$ given by the consensus algorithm~(\ref{eq:periodicConsensnsProt}).
Here, $r_{\asym}^T$ with $T=1$ does not necessarily correspond to the static-optimal solution in (\ref{eq:asymConvFact}).
\end{remark}

\section{Performance evaluations}\label{sec:numericalSimulations}

In this section, we demonstrate the effectiveness of the proposed method based on various numerical experiments. The number of datasets per learning was set to \num{1000}, and we performed the online learning (i.e., we set the batch size one). The weights of the edges were initialized to $0.1$. For evaluation, \num{10000} samples were used. Based on this setup, experiments were conducted in PyTorch \cite{pytorch} using Adam with a learning rate of 0.01 for training. Throughout our numerical experiments, we fixed the length of the time window for optimization to $T=10$ and adopted the periodic continuation introduced in Subsection~\ref{sec:periodic}.

To assess the effectiveness of the proposed approach, we compared the performance of the proposed method to that of two baseline strategies; the static optimal strategy in \cite{Xiao2004} and the finite-time distributed algorithm in \cite{Safavi2015}. These three approaches were compared on both empirical networks and random synthetic networks.

\subsection{Baseline strategies}

Here, we introduce the static optimal strategy presented in~\cite{Xiao2004} and the finite-time distributed algorithm in \cite{Safavi2015}.

\subsubsection{Static optimal strategy}

The static optimal strategy uses a set of time-invariant weights. Thus, the problem is to find such a set of weights that together minimize the asymptotic convergence factor.

Assume that an initial state $x_0$ is a deterministic vector. Additionally, assume that the edge weights $w_{ij}(k)$ do not depend on the time~$k$. Under these assumptions, the authors of~\cite{Xiao2004} showed that the problem of finding the static edge weights minimizing the (worst-case) asymptotic convergence factor~\eqref{eq:asymConvFact} can be reduced to solving a linear matrix inequality, which can be solved globally and efficiently \cite{Boyd1994}. As a baseline strategy, we adopted the
following time-invariant consensus protocol:
\begin{equation}\label{eq:lti}
  \begin{multlined}
    x_i(k+1) = x_i(k) + \sum_{j\in \mathcal N_i} w^{\stat}_{ij}(x_j(k) - x_i(k)),
    \ \
    0\leq k\leq T-1,
  \end{multlined}
\end{equation}
where $w^{\stat}_{ij}$ are the static optimal weights obtained by solving the linear matrix inequality.

\subsubsection{Finite-time distributed algorithm}

The finite-time distributed algorithm uses a method called successive nulling eigenvalues. Basically, it computes diagonal weights for a given set of non-diagonal weights so as to achieve a finite-time consensus.
Let $W = [w_{ij}]_{i, j}$ be a symmetric $N\times N$ matrix
satisfying $w_{ij}=0$ for all $(i,j) \notin E$.
Let $\lambda_1$, \dots, $\lambda_K$ be the distinct eigenvalues of $W$.
Then, the finite-time distributed algorithm proposed in \cite{Safavi2015} can be written as
\begin{equation}\label{eq:ftd}
  \begin{multlined}
    x_i(k+1) = a(k)x_i(k) + \sum_{j\in \mathcal N_i} w_{ij}x_j(k),
    \ \
    0\leq k\leq K-1,
  \end{multlined}
\end{equation}
where
$a(k) = \lambda_{k+1}$ for $k = 0, ..., K-2$ and
\begin{align}\label{eq:aK-1}
  a(K-1) = \dfrac{1}{(\lambda_{K-1}-\lambda_K)\cdots(\lambda_1-\lambda_K)}+\lambda_K.
\end{align}
The authors of~\cite{Safavi2015} showed that, under the assumption that the eigenvector corresponding to the eigenvalue $\lambda_K$ is $\onev$, the finite-time distributed algorithm~\eqref{eq:ftd} achieves an average consensus at a time~$K$ (i.e., $\epsilon_K = 0$). To satisfy the assumption on the matrix~$W$, we let $W$ be the Laplacian matrix of the graph~$G$ in this numerical experiment.

\subsection{Empirical networks}\label{subsec:deterministic}

\begin{table}[tb]
  \sisetup{
    round-mode = places,
    round-precision = 1,
    detect-all = true,
    detect-inline-weight = math}%  
  \begin{center}
    \caption{Numbers of nodes ($N$), edges ($M$), and distinct eigenvalues of the Laplacian matrices ($K$) for the empirical networks}
    \label{table:empiricalNets}
    \begin{tabular}{lccc}
      \toprule
                                                                              &
      $N$                                                                     &
      $M$                                                                     &
      $K$
      \\
      \midrule
      {Krackhardt kite}                                                       &
      10 &
      18 &
      10
      \\
      {Chvat\'al}                                                             &
      12         &
      24         &
       7
      \\
      {Pappus}                                                                &
      18          &
      27          &
       5
      \\
      {Davis}                                                                 &
      32           &
      89           &
      32
      \\
      {Karate}                                                                &
      34          &
      78          &
      30
      \\
      {Tutte}                                                                 &
      46           &
      69           &
      31
      \\
      \bottomrule
    \end{tabular}
  \end{center}
\end{table}

In this subsection, we consider the following six empirical and synthetic deterministic networks: Krackhardt kite graph (\textit{Krackhardt kite}), Chv\'atal graph (\textit{Chv\'atal}), Puppus graph (\textit{Puppus}), bipartite network of Southern women and clubs (\textit{Davis}), social network of a Karate club~\cite{Zachary1977} (\textit{Karate}), and Tutte graph (\textit{Tutte}).
We summarize the numbers of nodes, edges, and distinct eigenvalues of the Laplacian matrices of the empirical networks in Table~\ref{table:empiricalNets}.

\begin{figure}
  \centering
  \includegraphics[width=.8\linewidth]{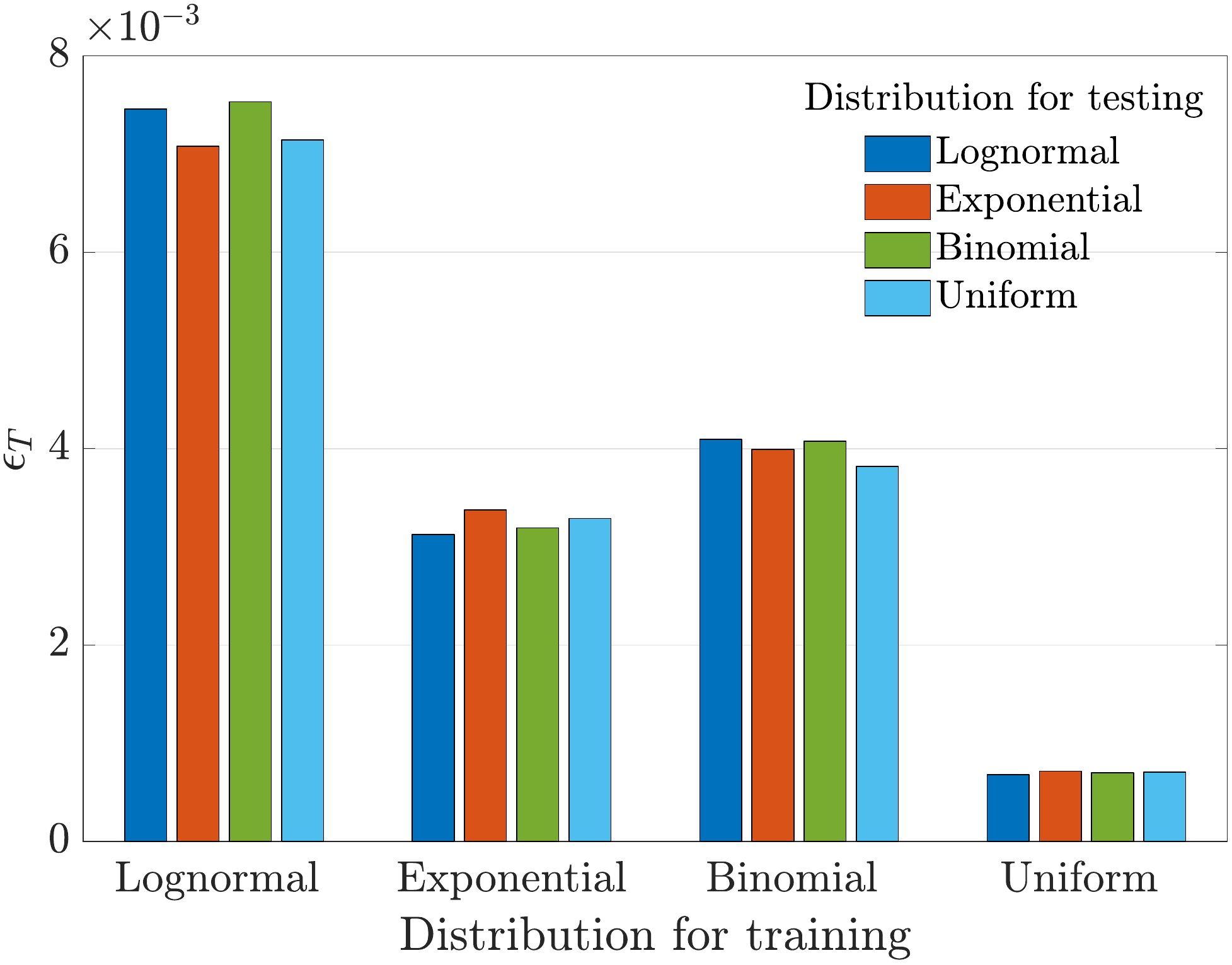}
  \caption{Performance $\varepsilon_T$ of the Karate network trained using four different distributions (lognormal, exponential, binomial, and uniform). The performance value achieved using the baseline strategy are between 0.2 and 0.22 for all distributions.}
  \label{fig:distributions}
\end{figure}

\subsubsection{Effects of initial distributions}

First, we trained a Karate network with $T=10$ based on four types of initial states sampled from 1) log-normal distribution with parameters $\mu = 0$ and $\sigma=1.5$, 2) exponential distribution with mean~$1$, 3) binomial distribution with parameters $n=50$ and $p=0.5$, and 4) uniform distribution. We then evaluated the performance $\varepsilon_{\,T}$ for each of the trained temporal networks by using the aforementioned distributions to generate the initial state $x(0)$. We present the results of this experiment in Fig.~\ref{fig:distributions}. One can see that the networks trained using different distributions perform differently. Although performance is robust to the choice of an initial distribution overall, the network trained using a uniform distribution yields the best results.
  Therefore, from now on, we will assume that the initial state of each node independently follows a uniform distribution on the interval~$[-1, 1]$. 

\subsubsection{Dynamically changing edge weights}

\begin{figure}[tb]
  \newcommand{\networkWidth}{.23\linewidth}
  \centering
  \begin{subfigure}{\networkWidth}
    \centering
    \includegraphics[width=.8\linewidth]{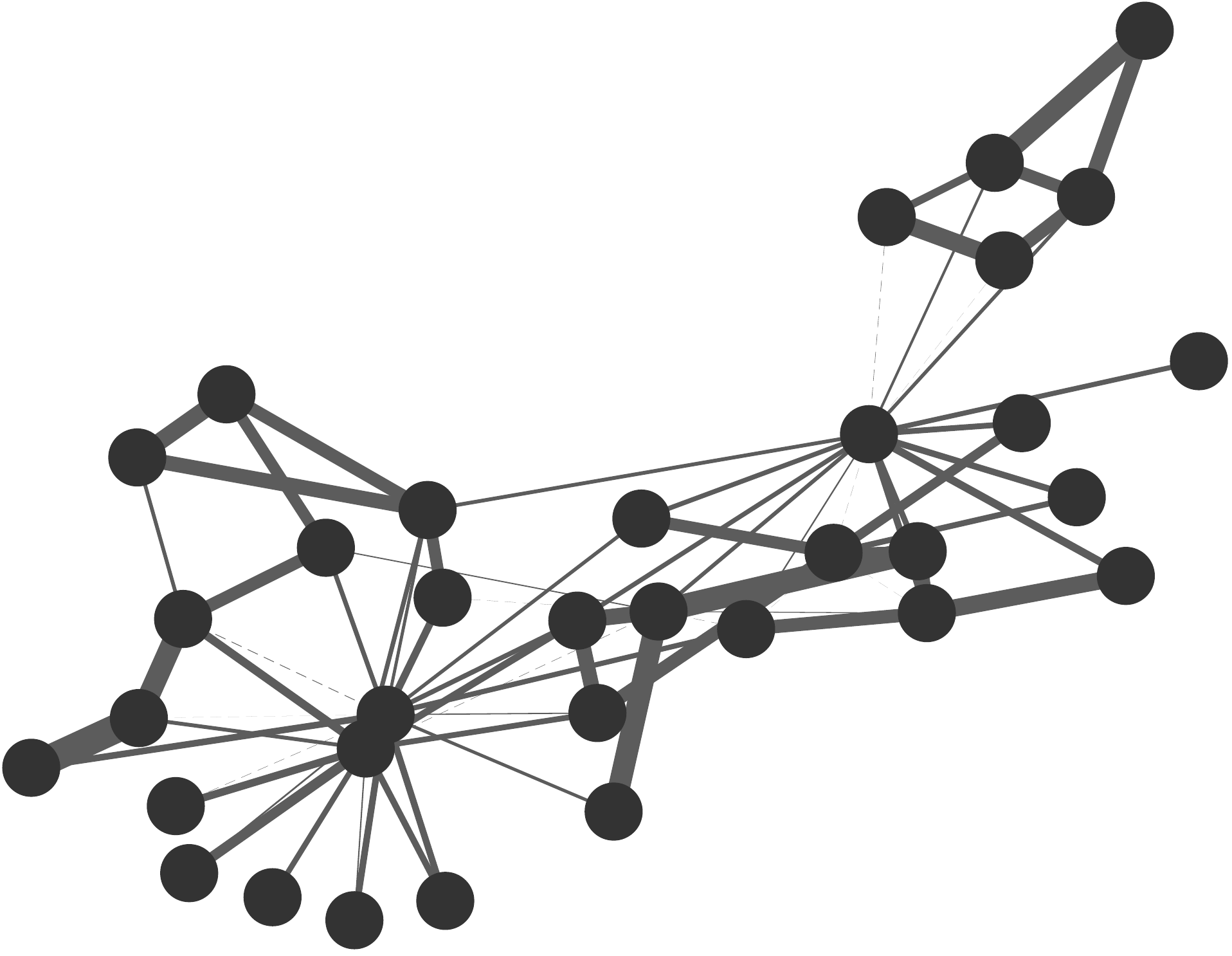}
    \caption{$k=0$}
  \end{subfigure}
  \hfil
  \begin{subfigure}{\networkWidth}
    \centering
    \includegraphics[width=.8\linewidth]{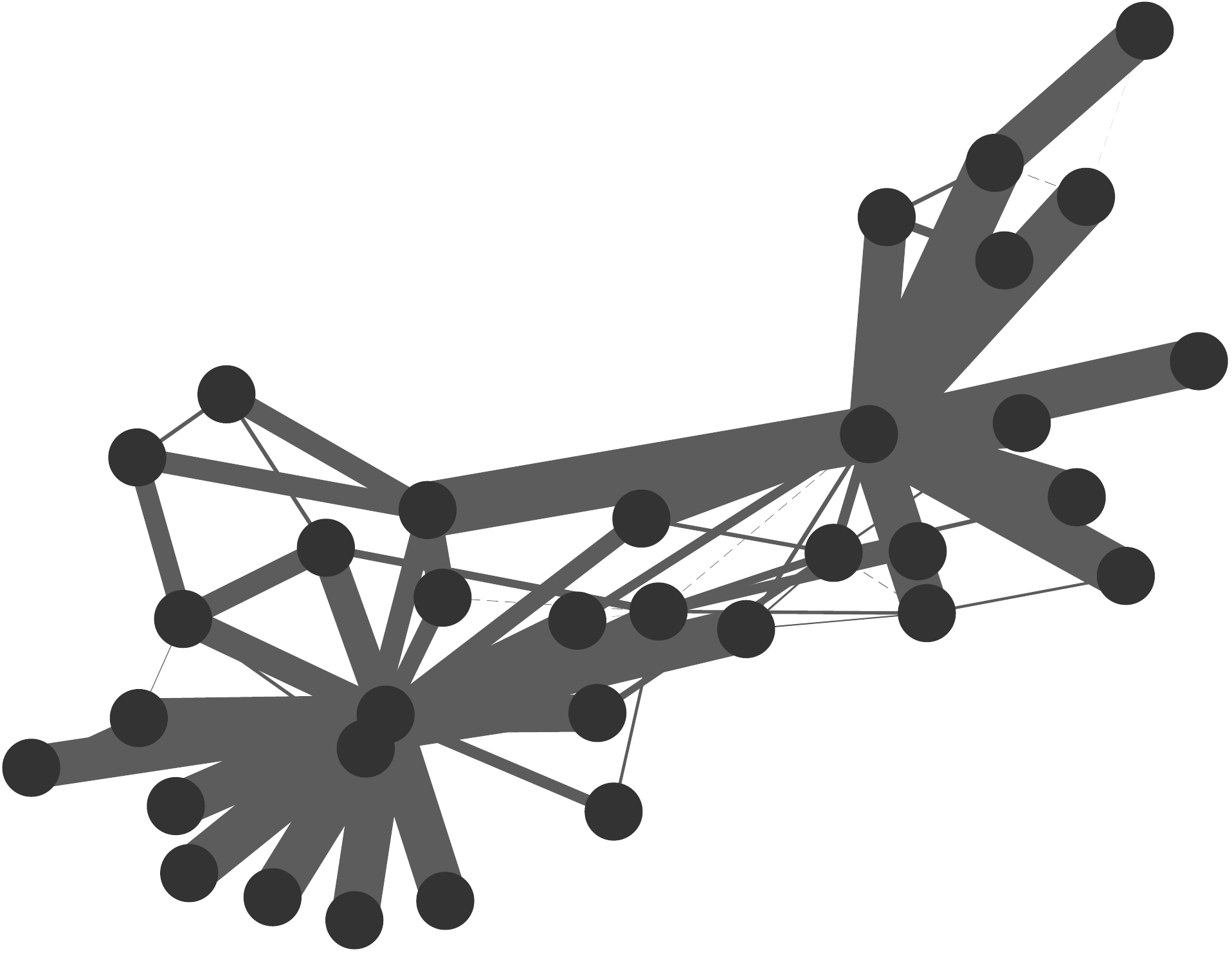}
    \caption{$k=1$}
  \end{subfigure}
  \hfil
  \begin{subfigure}{\networkWidth}
    \centering
    \includegraphics[width=.8\linewidth]{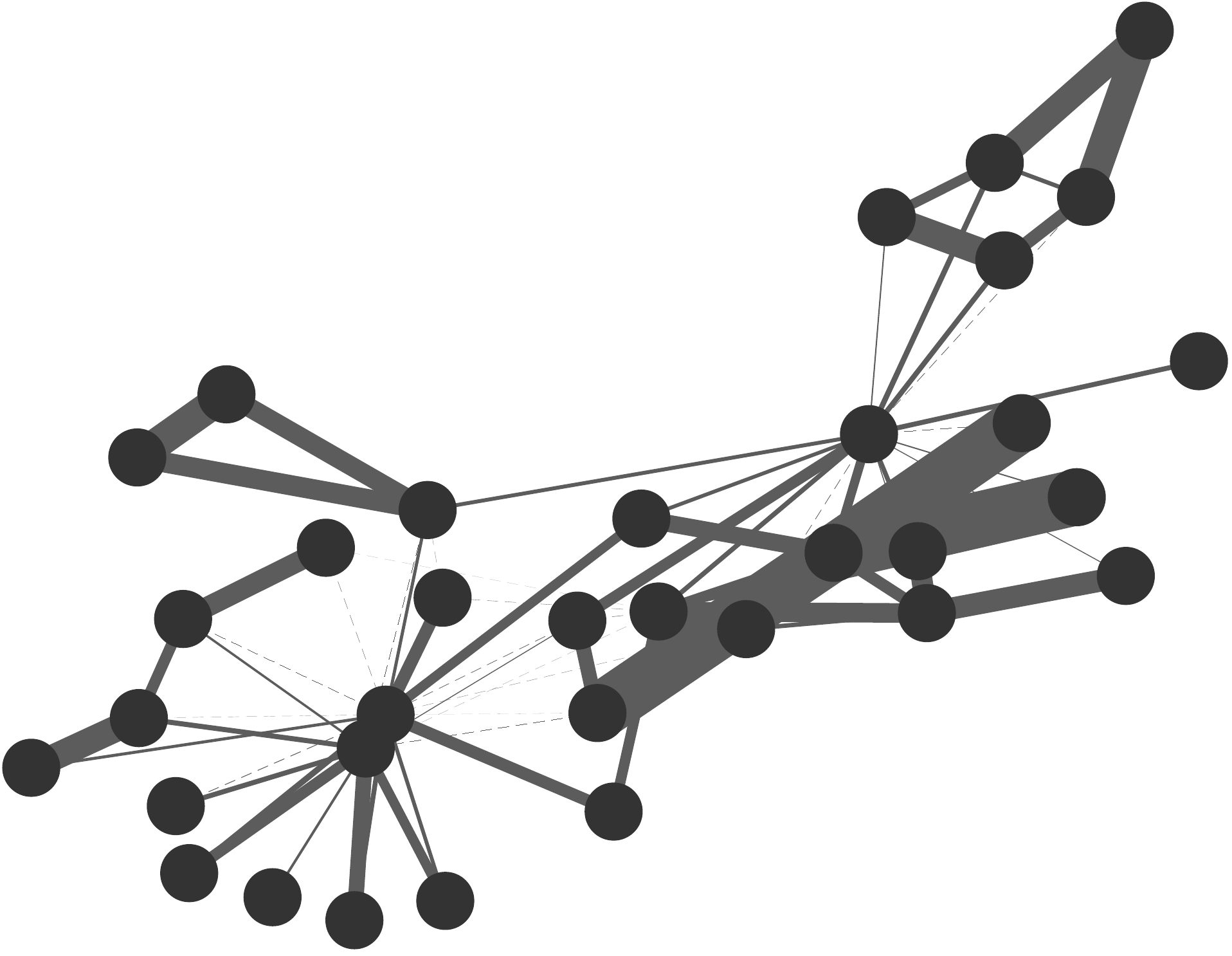}
    \caption{$k=2$}
  \end{subfigure}
  \hfil
  \begin{subfigure}{\networkWidth}
    \centering
    \includegraphics[width=.8\linewidth]{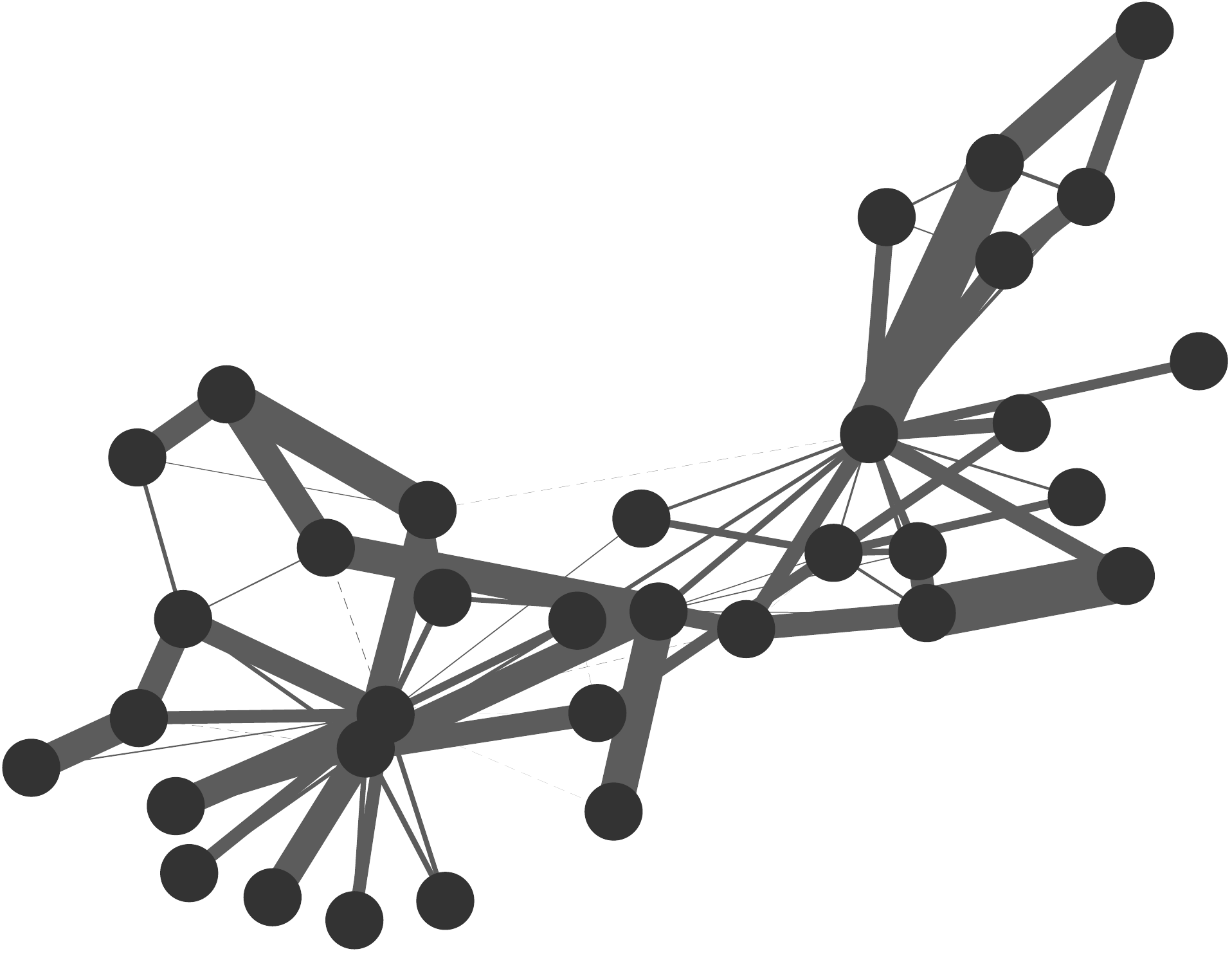}
    \caption{$k=3$}
  \end{subfigure}
  \\
  \vspace{2mm}
  \begin{subfigure}{\networkWidth}
    \centering
    \includegraphics[width=.8\linewidth]{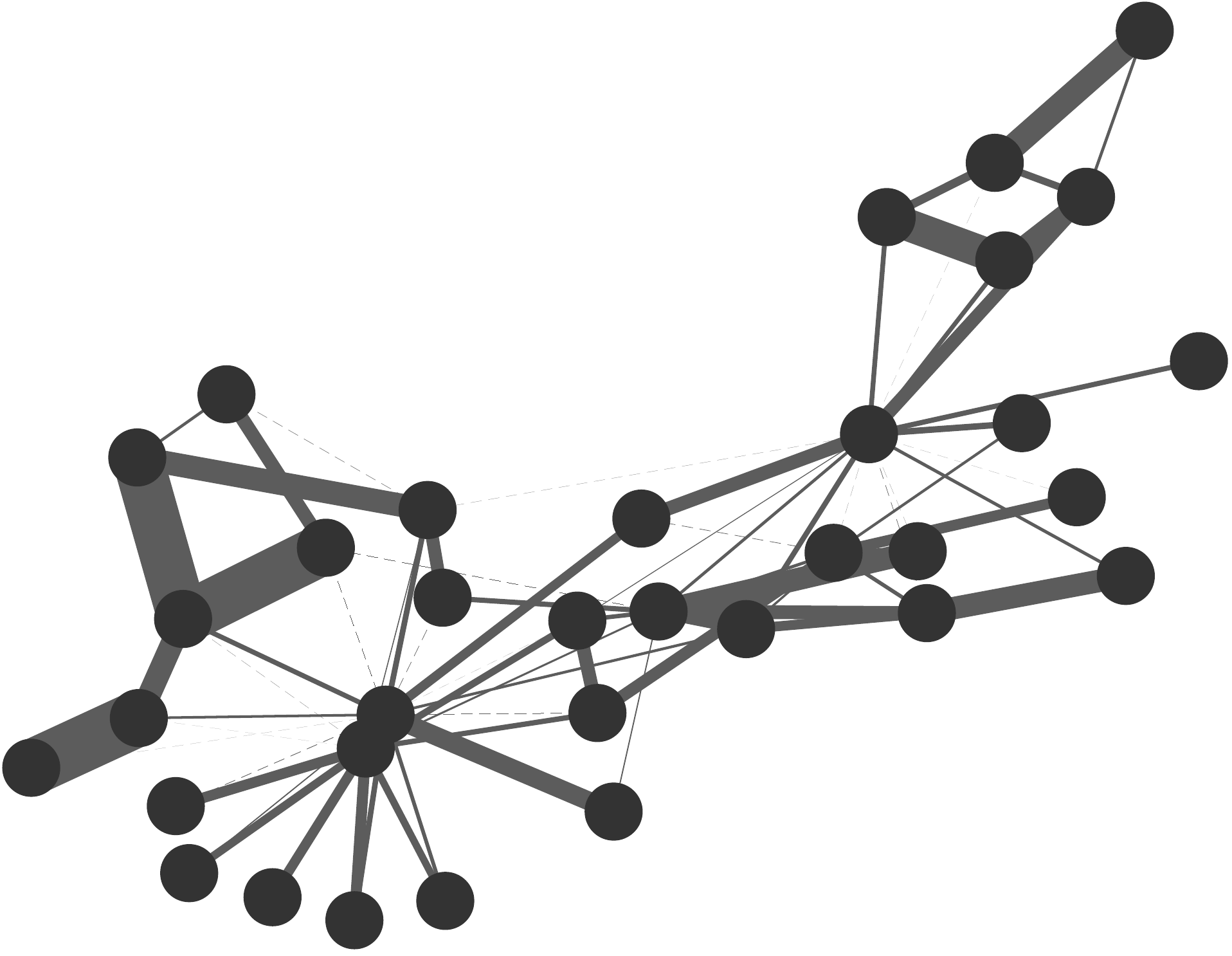}
    \caption{$k=4$}
  \end{subfigure}
  \hfil
  \begin{subfigure}{\networkWidth}
    \centering
    \includegraphics[width=.8\linewidth]{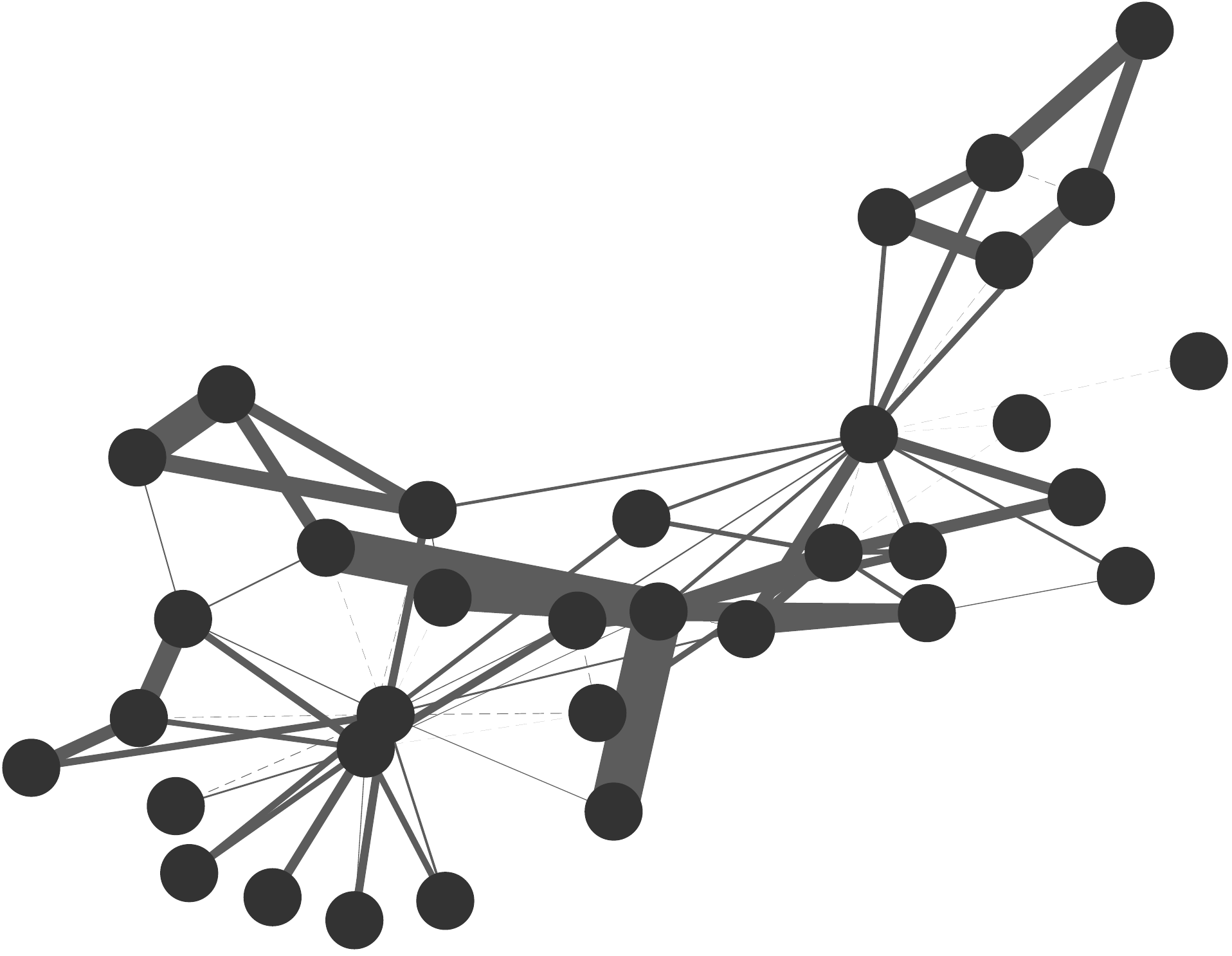}
    \caption{$k=5$}
  \end{subfigure}
  \hfil
  \begin{subfigure}{\networkWidth}
    \centering
    \includegraphics[width=.8\linewidth]{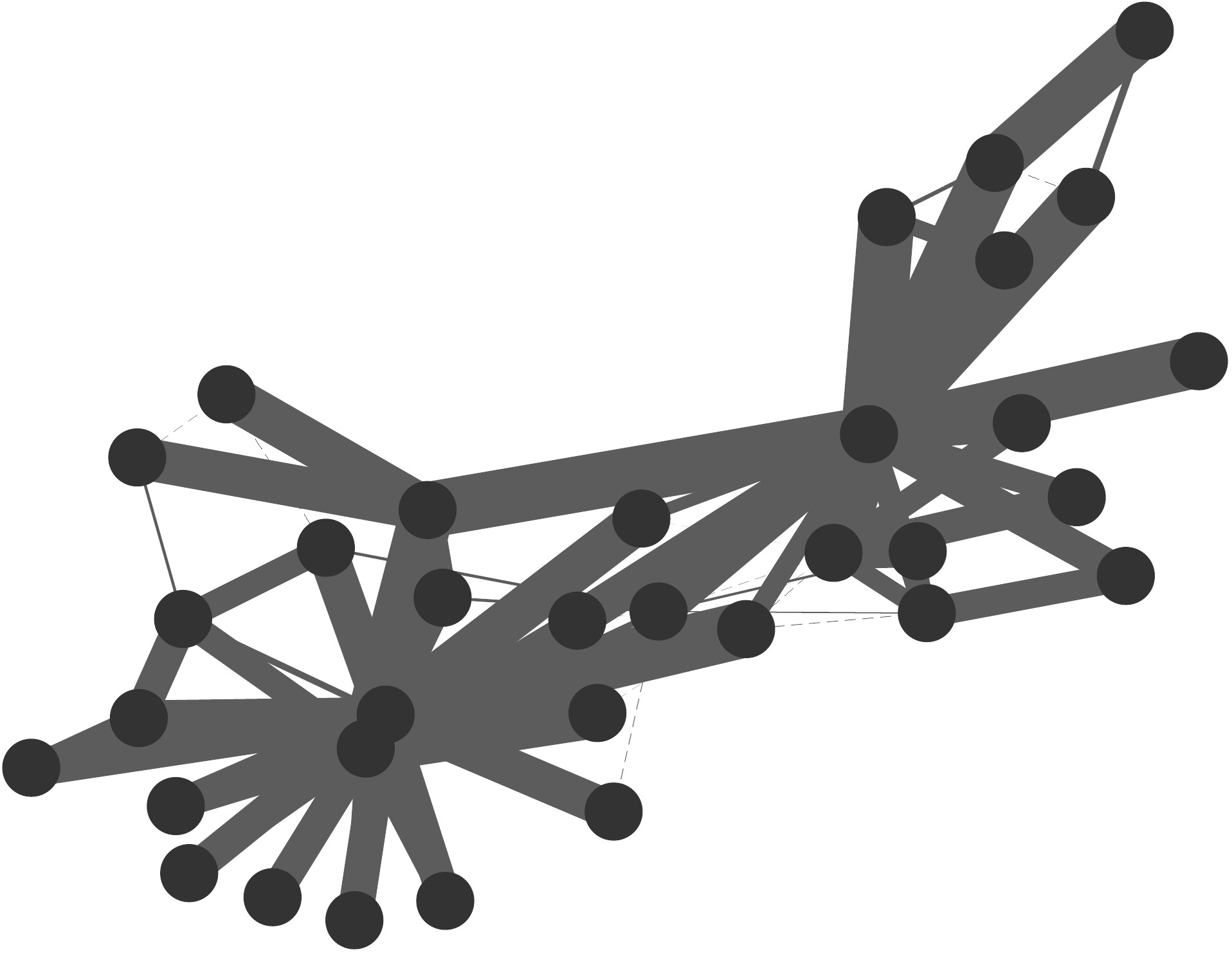}
    \caption{$k=6$}
  \end{subfigure}
  \hfil
  \begin{subfigure}{\networkWidth}
    \centering
    \includegraphics[width=.8\linewidth]{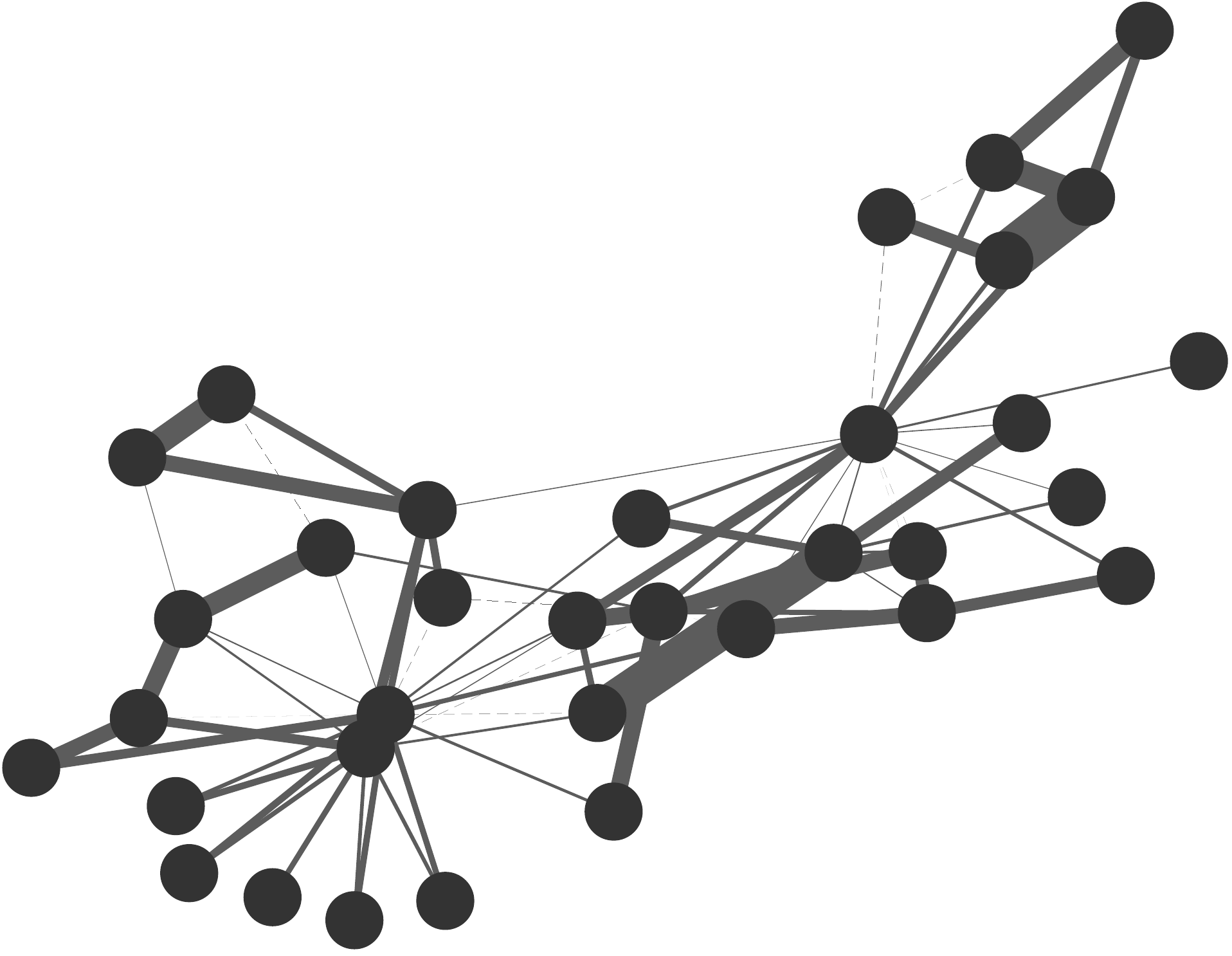}
    \caption{$k=7$}
  \end{subfigure}
  \\
  \vspace{2mm}
  \begin{subfigure}{\networkWidth}
    \centering
    \includegraphics[width=.8\linewidth]{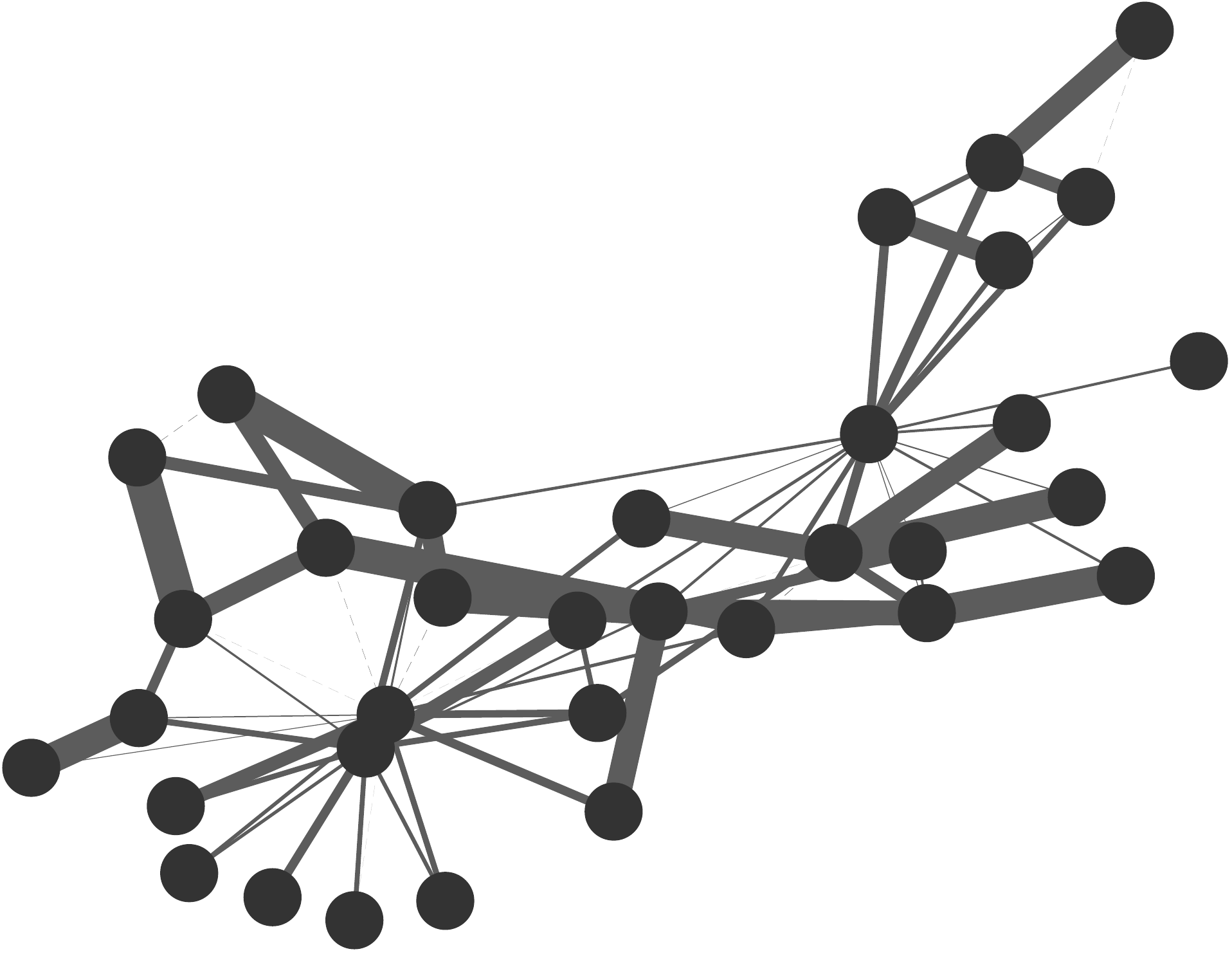}
    \caption{$k=8$}
  \end{subfigure}
  \hfil
  \begin{subfigure}{\networkWidth}
    \centering
    \includegraphics[width=.8\linewidth]{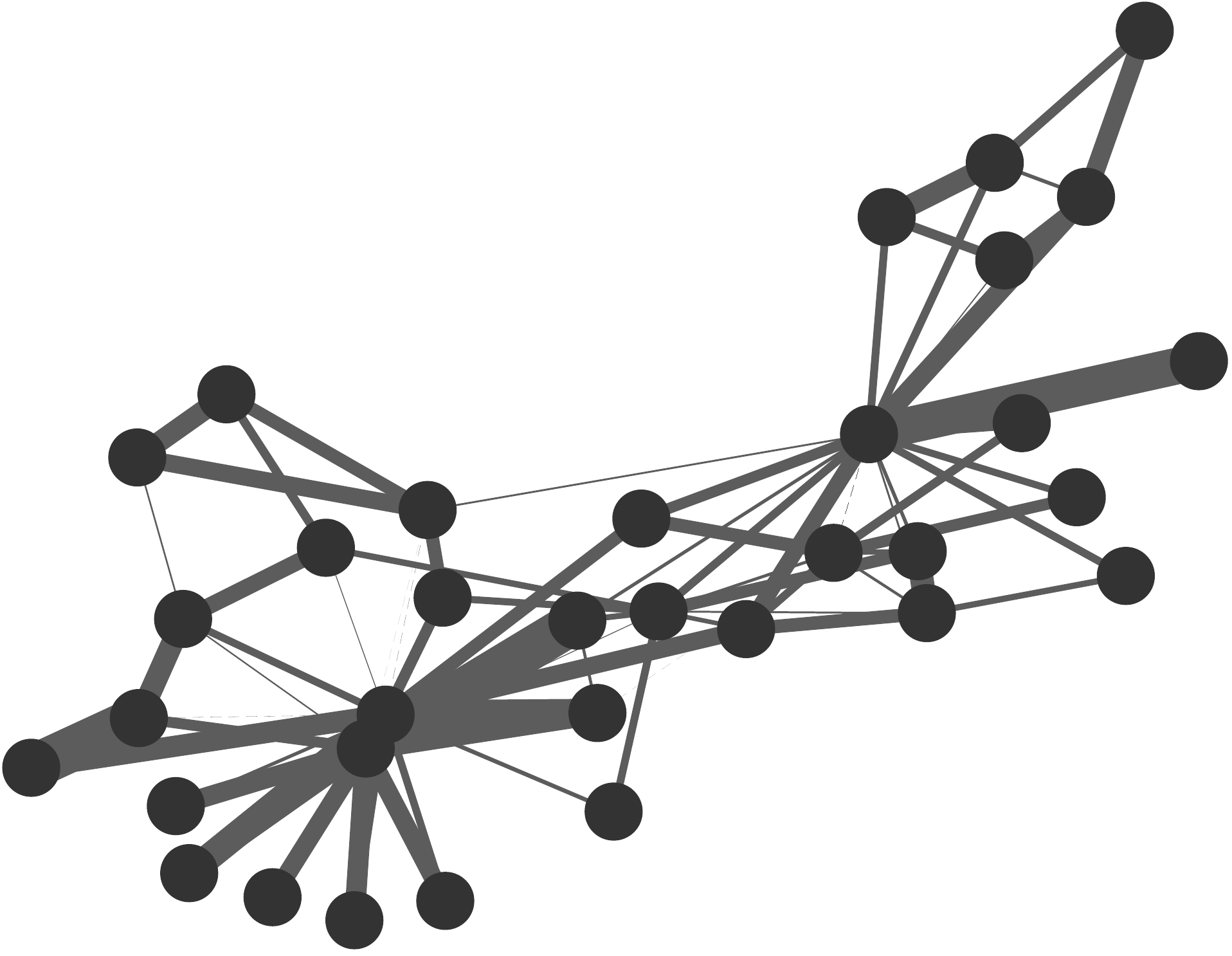}
    \caption{$k=9$}
  \end{subfigure}
  \hfil
  \begin{subfigure}{.46\linewidth}
    \centering
    \includegraphics[width=.4\linewidth]{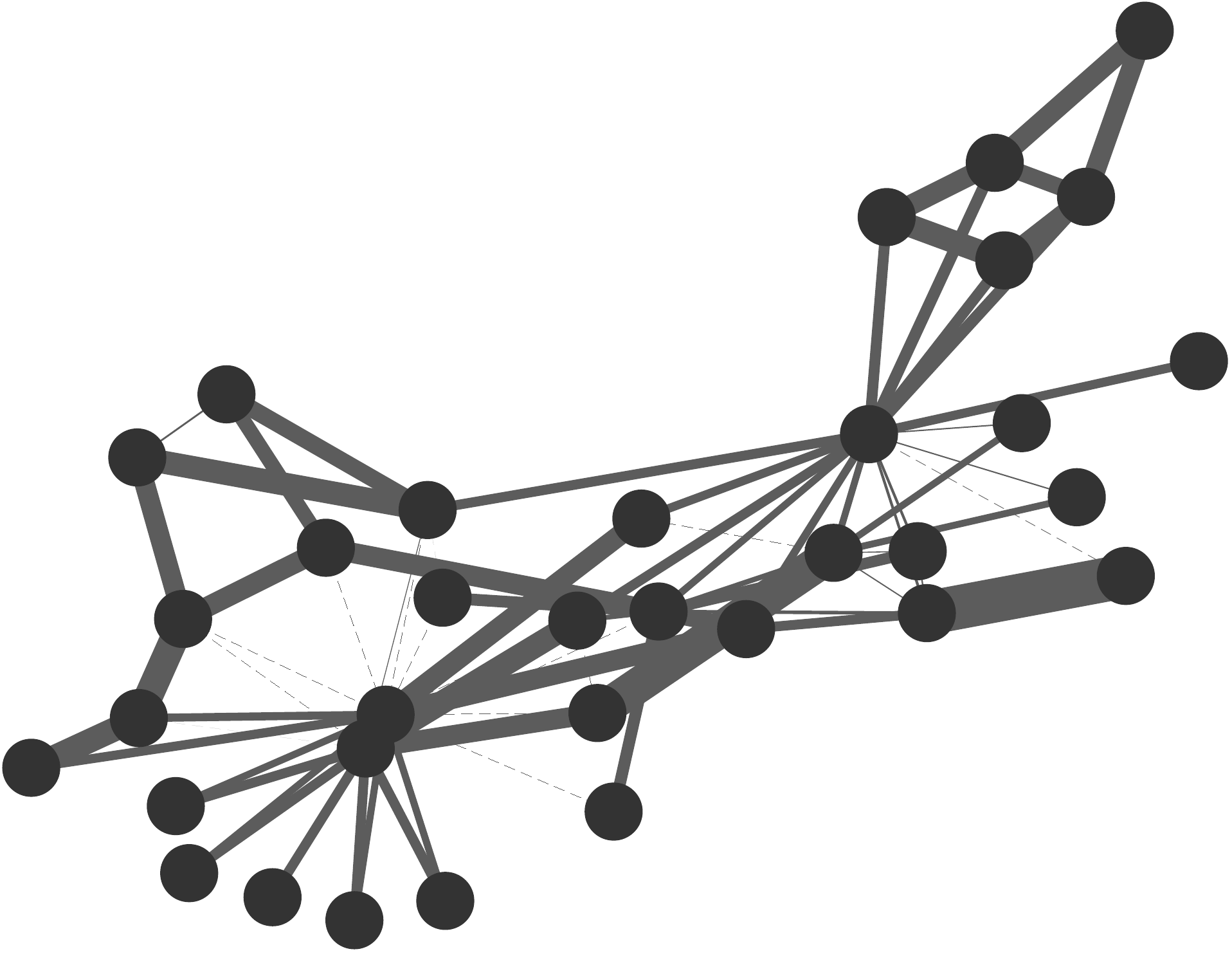}
    \caption{Static optimal}
  \end{subfigure}
  \caption{Optimized weights of edges for the Karate network. (a)--(j): Proposed method. (k): Static optimal strategy. The width of the lines indicates the values of the weights (thicker line = larger weight). Edges with weights less than $10^{-2}$ are indicated by dashed lines.}
    \label{fig:weights:karate}
\end{figure}

Before comparing the proposed approach to the baseline approaches, we show the optimized edge weights for the Karate network
  at times $k=0, \dotsc, 9$ by minimizing $\epsilon_T$ with $T=10$ in Fig.~\ref{fig:weights:karate}.
One can see that the optimized edge weights of the network change dynamically in a nontrivial manner.

\begin{table}[tb]
  \sisetup{
    round-mode = places,
    round-precision = 1,
    detect-all = true,
    detect-inline-weight = math}%  
  \begin{center}
    \caption{Consensus errors~$\varepsilon_K$ for empirical networks}
    \label{table:}
    \begin{tabular}{lccc}
      \toprule
                                                                                            &
      Proposed                                                                              &
      \begin{tabular}{@{}c@{}}Finite-time consensus\end{tabular}                                                            &
      \begin{tabular}{@{}c@{}}Static optimal\end{tabular}                                                              \\
      \midrule
      {Krackhardt kite}                                                                     &
      \num{1.4625e-05}              &
      \textbf{\num{6.3620e-10}} &
      \num{1.6569e-01}
      \\
      {Chvat\'al}
                                                                                            &
      \num{2.1052e-05}                        &
      \textbf{\num{3.6310e-12}}           &
      \num{3.3397e-03}
      \\
      {Pappus}                                                                              &
      \num{1.0018e-03}                       &
      \textbf{\num{1.8452e-13}}          &
      \num{1.7366e-01}
      \\
      {Davis}                                                                               &
      \textbf{\num{1.5292e-07}}                 &
      \num{3.1629e+17}                      &
      \num{3.9184e-03}
      \\
      {Karate}                                                                              &
      \textbf{\num{1.8899e-07}}                &
      \num{1.2924e+18}                     &
      \num{1.2096e-01}
      \\
      {Tutte}                                                                               &
      \textbf{\num{2.0027e-05}}               &
      \num{5.6992e+02}                    &
      \num{2.5061e-01}
      \\
      \bottomrule
    \end{tabular}
  \end{center}
\end{table}

\subsubsection{Consensus errors for a finite-time window of length $K$}

We empirically evaluated the average consensus error~$\epsilon_K$, where $K$ denotes the numbers of distinct eigenvalues of the Laplacian matrices of the networks. We also applied the static optimal strategy and finite-time consensus algorithm and empirically evaluated their average consensus errors $\epsilon_K$. We list these consensus errors in Table~\ref{table:}. The proposed method achieves small consensus errors and outperforms the static optimal strategy for all of the networks. The finite-time consensus algorithm achieves accurate consensuses for small networks ({Krackhardt kite}, {Chv\'atal}, and {Puppus}). However, this algorithm fails to achieve a consensus for large networks. This can be attributed to numerical instability when computing the coefficient~{\eqref{eq:aK-1}}, whose magnitude tends to become significantly smaller as the network size (and number $K$ of distinct eigenvalues of the Laplacian matrix) increases. For computing the coefficient~\eqref{eq:aK-1}, we adopted the following straightforward procedure. We first numerically computed the eigenvalues~$\lambda_1$, \dots, $\lambda_K$. We then computed the product of the differences~$\lambda_{K-1}-\lambda_K$, \dots, $\lambda_1-\lambda_K$. Finally, we calculated the inverse of this product to obtain $a(K-1)$.

\begin{figure}[tb]
  \newcommand{\networkWidth}{.49\linewidth}
  \centering
  \begin{subfigure}{\networkWidth}
    \centering
    \includegraphics[width=\linewidth]{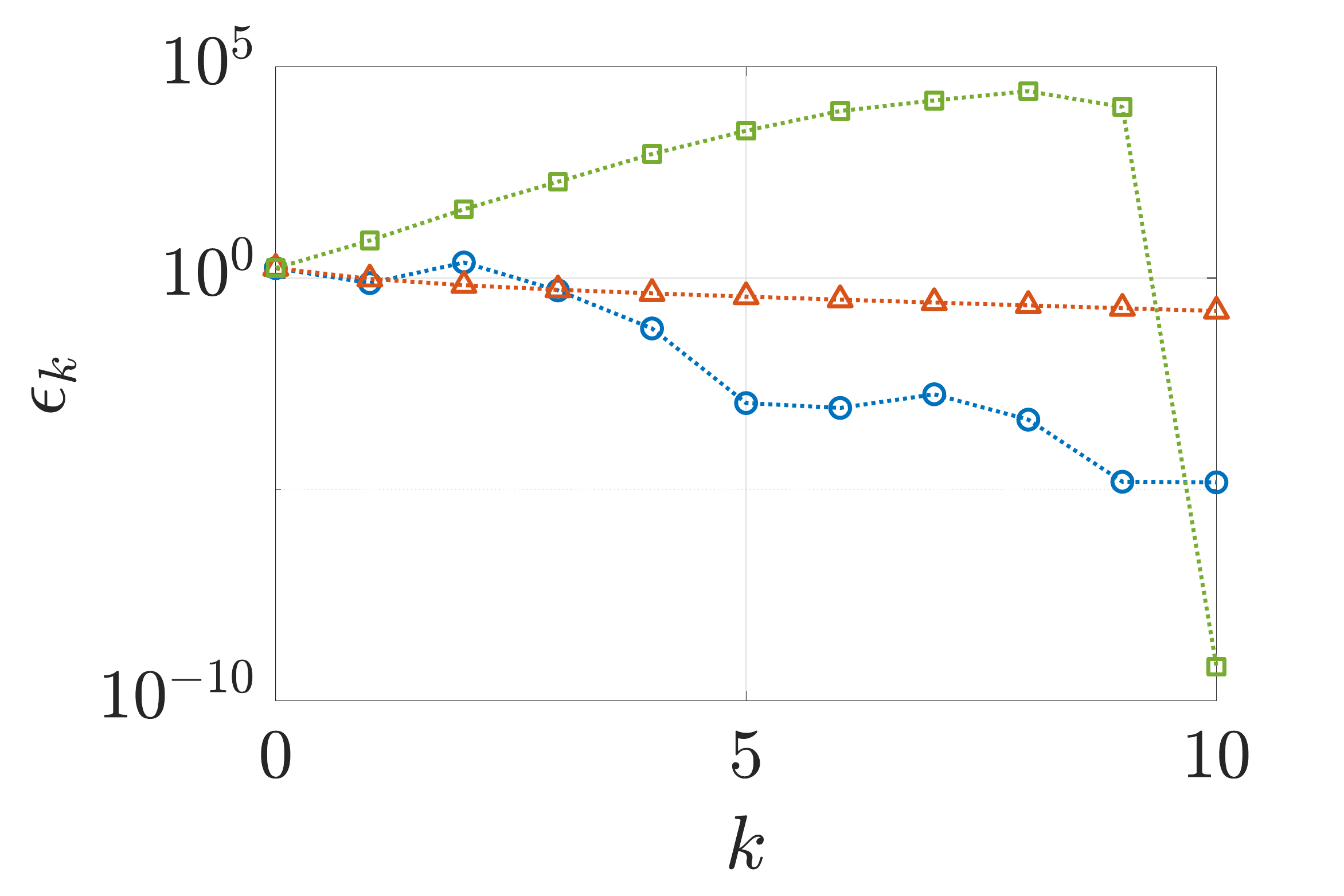}
    \caption{Krackhardt kite}
  \end{subfigure}
  \hfil
  \begin{subfigure}{\networkWidth}
    \centering
    \includegraphics[width=\linewidth]{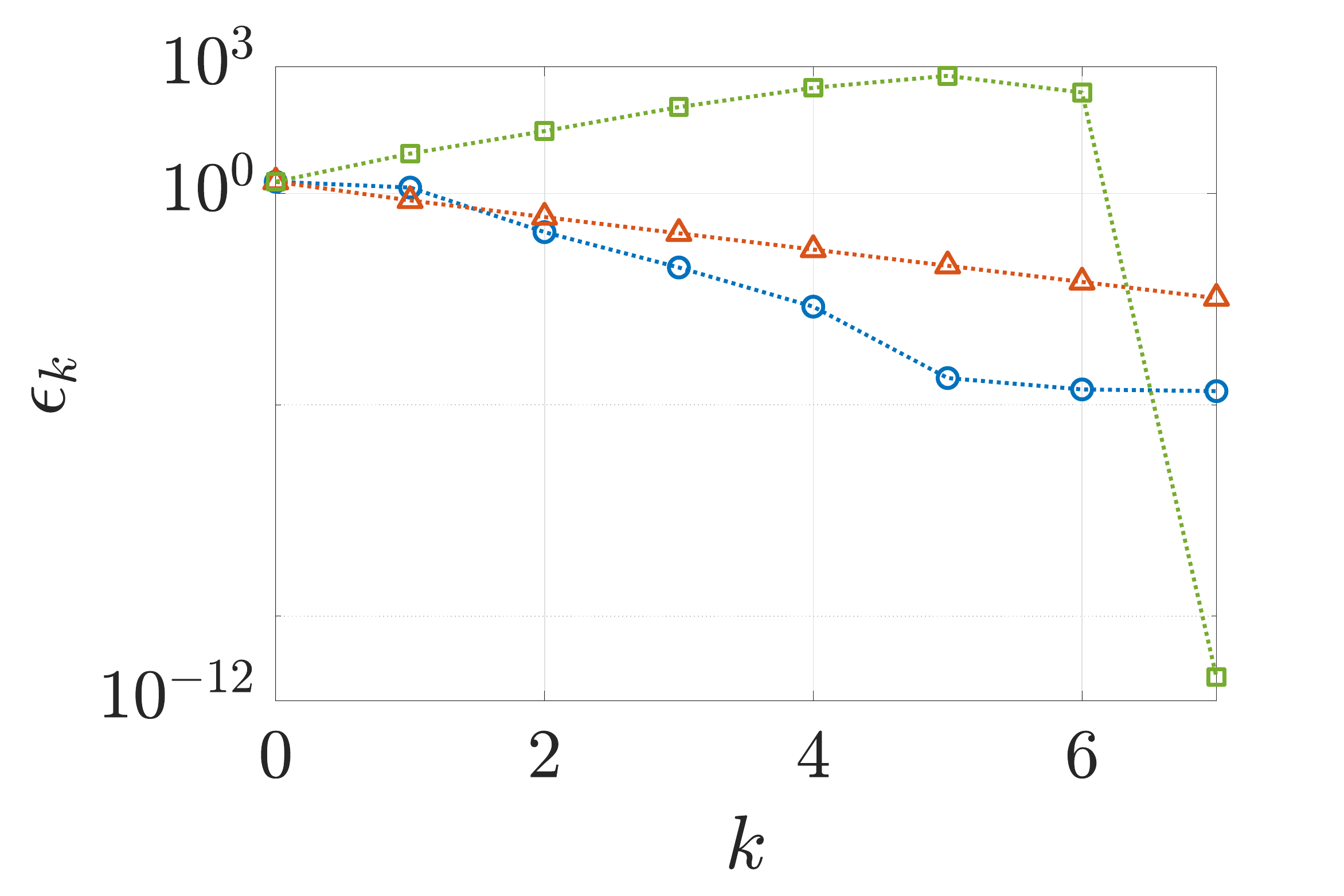}
    \caption{Chvatal}
  \end{subfigure}
  \\
  \vspace{2mm}
  \begin{subfigure}{\networkWidth}
    \centering
    \includegraphics[width=\linewidth]{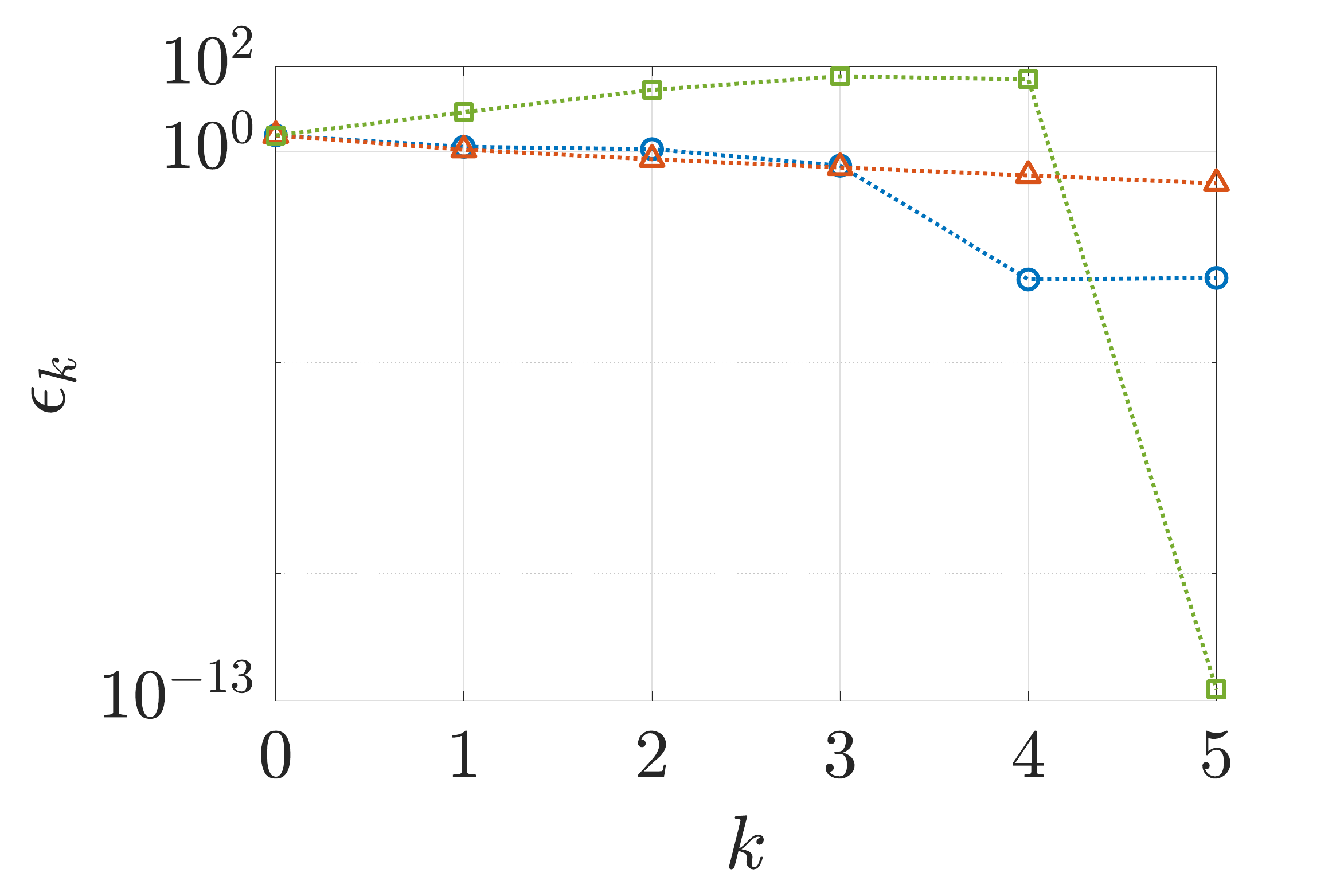}
    \caption{Pappus}
  \end{subfigure}
  \hfil
  \begin{subfigure}{\networkWidth}
    \centering
    \includegraphics[width=\linewidth]{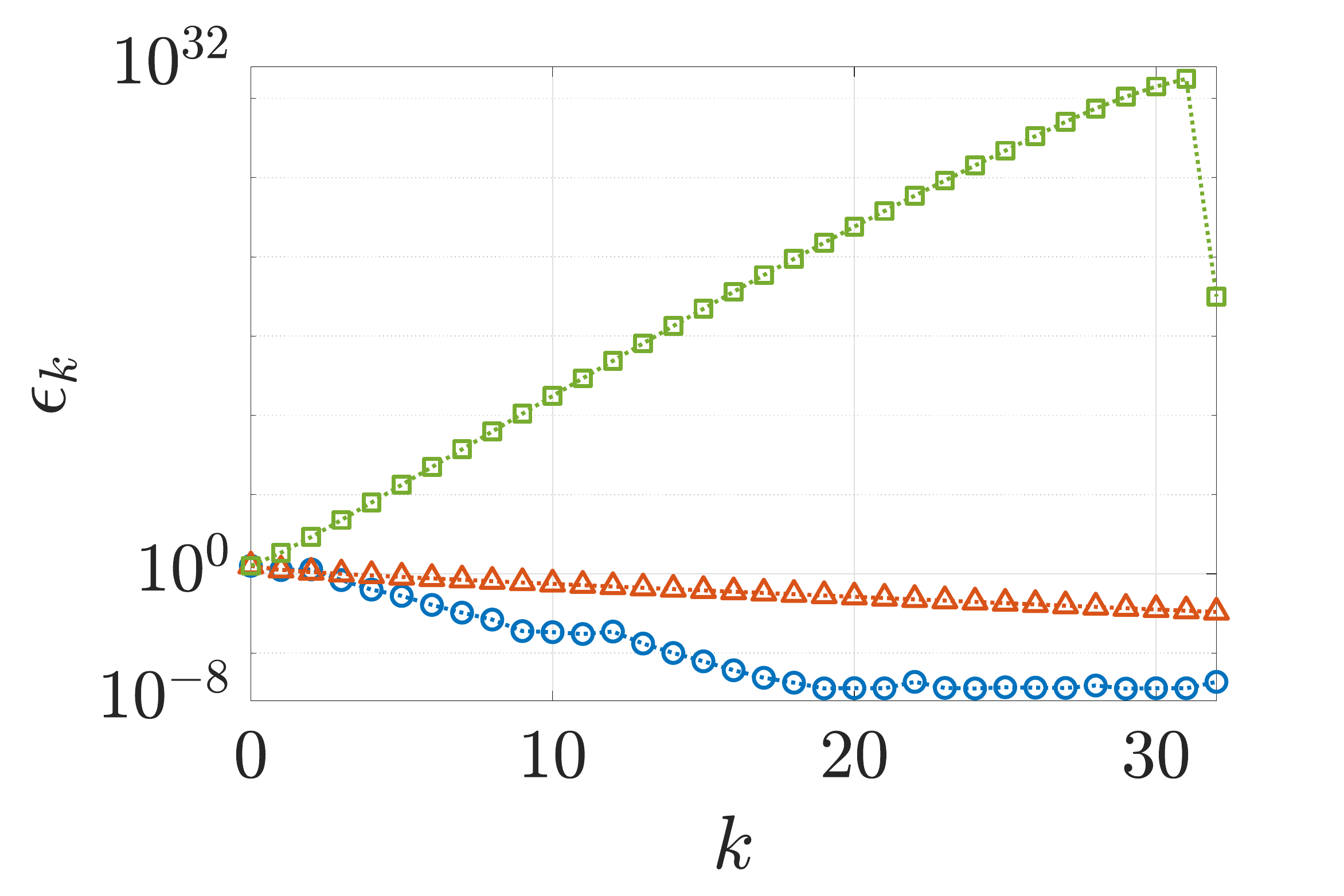}
    \caption{Davis}
  \end{subfigure}
  \\
  \vspace{2mm}
  \begin{subfigure}{\networkWidth}
    \centering
    \includegraphics[width=\linewidth]{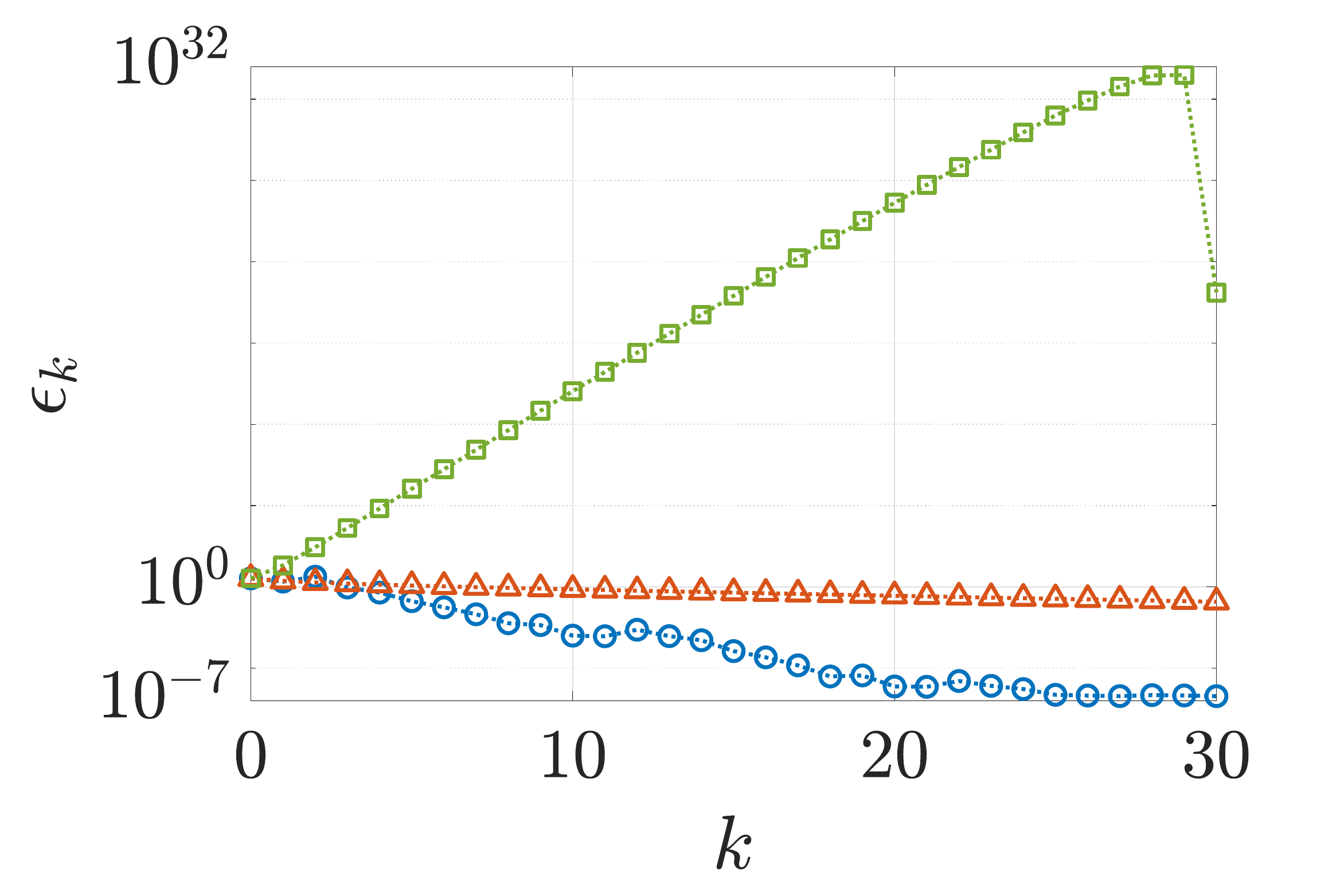}
    \caption{Karate}
  \end{subfigure}
  \hfil
  \begin{subfigure}{\networkWidth}
    \centering
    \includegraphics[width=\linewidth]{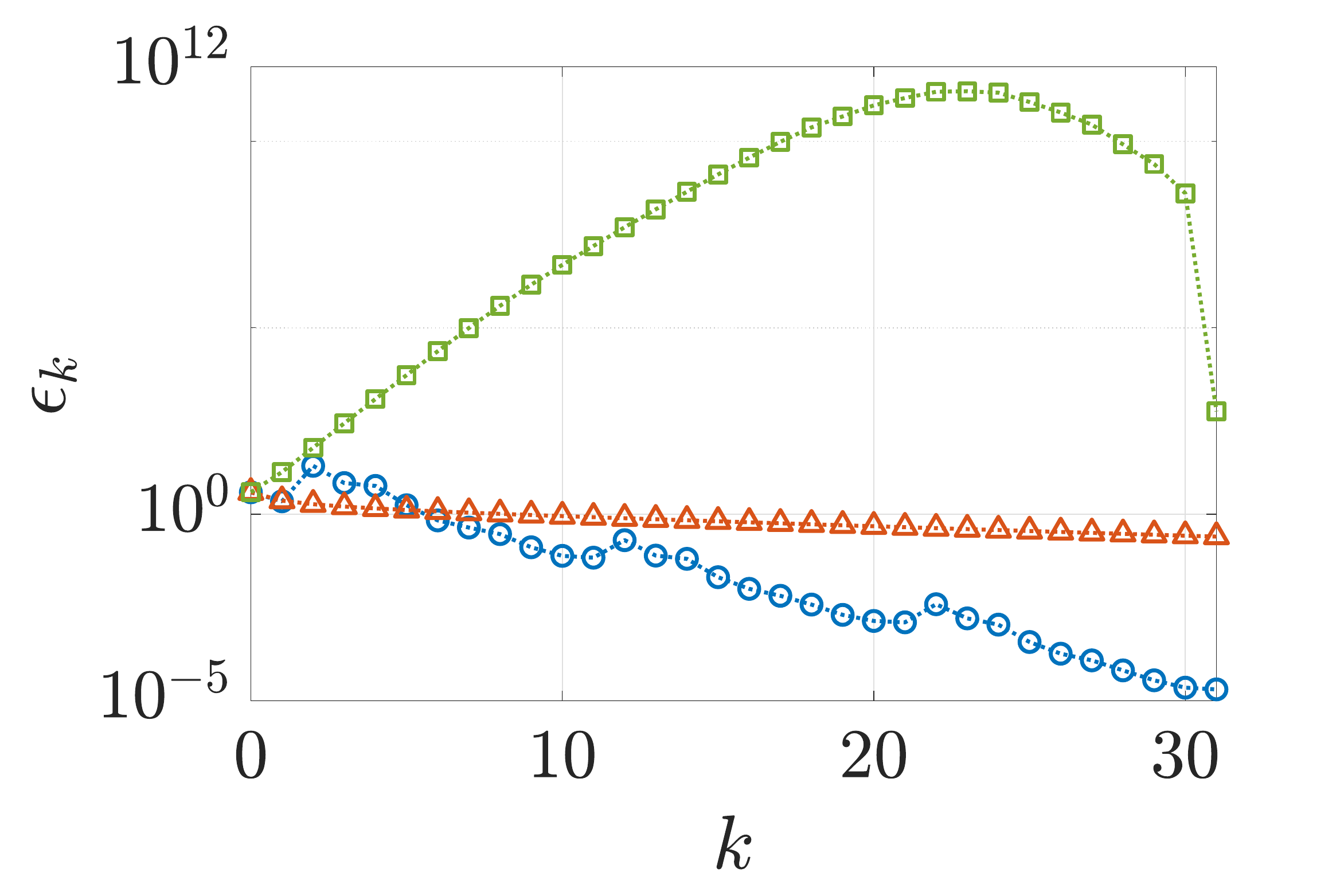}
    \caption{Tutte}
  \end{subfigure}
  \caption{Empirical averages of consensus errors. Circles: proposed method. Triangles: static optimal. Squares: finite-time consensus. All methods are terminated at time $k=K$, where the finite-time distributed algorithm is theoretically expected to achieve an exact average consensus.}
  \label{fig:consensus error trajectories empirical}
\end{figure}

In Fig.~\ref{fig:consensus error trajectories empirical}, we present the time evolutions of the consensus errors for the six empirical networks by the proposed method, the static optimal strategy, and the finite-time consensus algorithm. We terminate all three algorithms at time $k=K$, i.e., at the time when the finite-time distributed algorithm should theoretically supposed to achieve an exact average consensus. One can see huge consensus errors in the middle of the finite-time consensus algorithm, regardless of the underlying network. In contrast, the proposed method allows us to achieve an accurate consensus for any network size without overshooting.

\begin{table}[tb]
  \newcommand{\showFactors}[2]{
    #1 &
    \textbf{\input{"#2_DL_factor.txt"}}
    & \input{"#2_boyd_factor.txt"}}
  \sisetup{
    round-mode = places,
    round-precision = 1,
    detect-all = true,
    detect-inline-weight = math}
  \begin{center}
    \caption{Asymptotic convergence factors in empirical networks}
    \label{table:convergenceFactors}
    \begin{tabular}{lcc}
      \toprule
       & Proposed ($r_{\asym}^{10}$) & Static optimal ($r_{\asym}$) \\
      \midrule
      
    Krackhardt kite &
    \textbf{\num{3.1326e-01}}
    & \num{8.6088e-01}                          \\
      
    Chvatál &
    \textbf{\num{3.2788e-01}}
    & \num{4.1384e-01}                                          \\
      
    Pappus &
    \textbf{\num{3.3522e-01}}
    & \num{6.5108e-01}                                            \\
      
    Davis &
    \textbf{\num{3.9135e-01}}
    & \num{8.3294e-01}                                              \\
      
    Karate &
    \textbf{\num{4.9232e-01}}
    & \num{9.2503e-01}                                            \\
      
    Tutte &
    \textbf{\num{6.7873e-01}}
    & \num{9.4784e-01}                                              \\
      \bottomrule
    \end{tabular}
  \end{center}
\end{table}

\subsubsection{Asymptotic convergence factors}

Table~\ref{table:convergenceFactors} lists the convergence factors for the empirical networks obtained by the proposed method with periodic continuation and the static optimal strategy. We used Proposition~\ref{LEM:} to compute the convergence factors. The proposed method improves upon the static optimal strategy for all of the networks, which is consistent with our observations of consensus errors.

\begin{figure}[tb]
  \centering
  \begin{subfigure}{.9\linewidth}
    \centering
    \includegraphics[width=\linewidth]{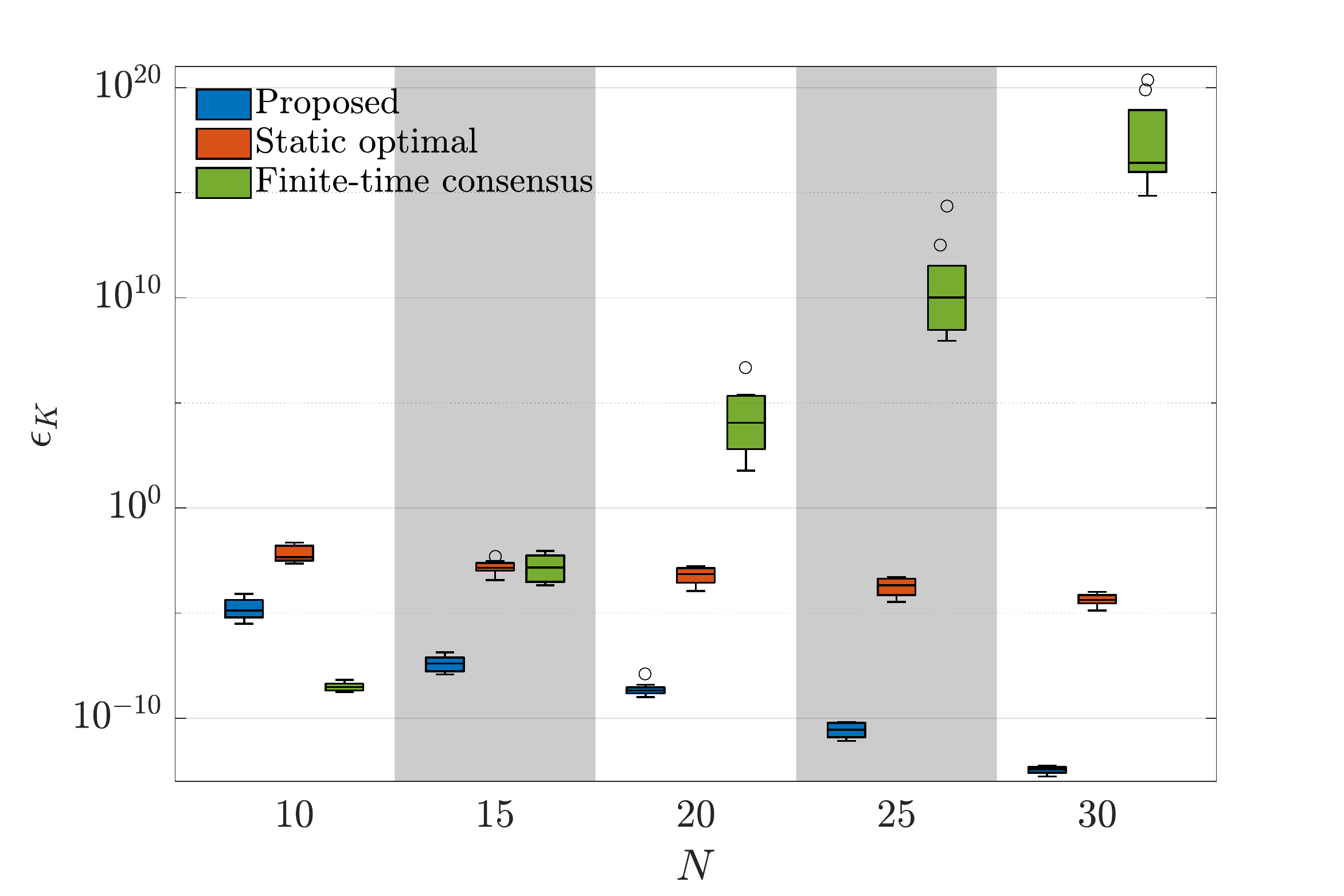}
    \caption{Barab\'asi-Albert network}
  \end{subfigure}
  \begin{subfigure}{.9\linewidth}
    \centering
    \includegraphics[width=\linewidth]{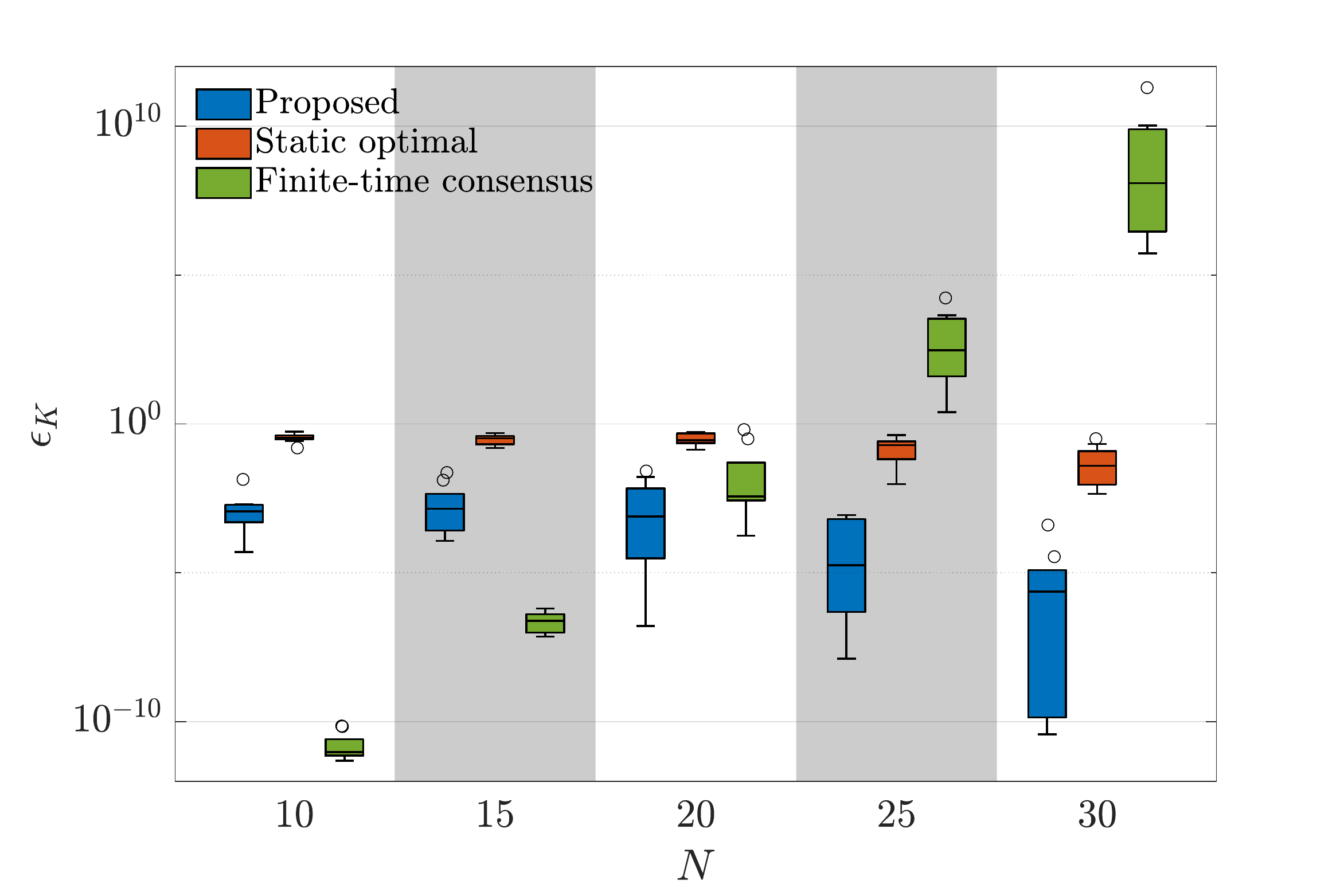}
    \caption{Erd\H{o}s-R\'enyi network}
  \end{subfigure}
  \begin{subfigure}{.9\linewidth}
    \centering
    \includegraphics[width=\linewidth]{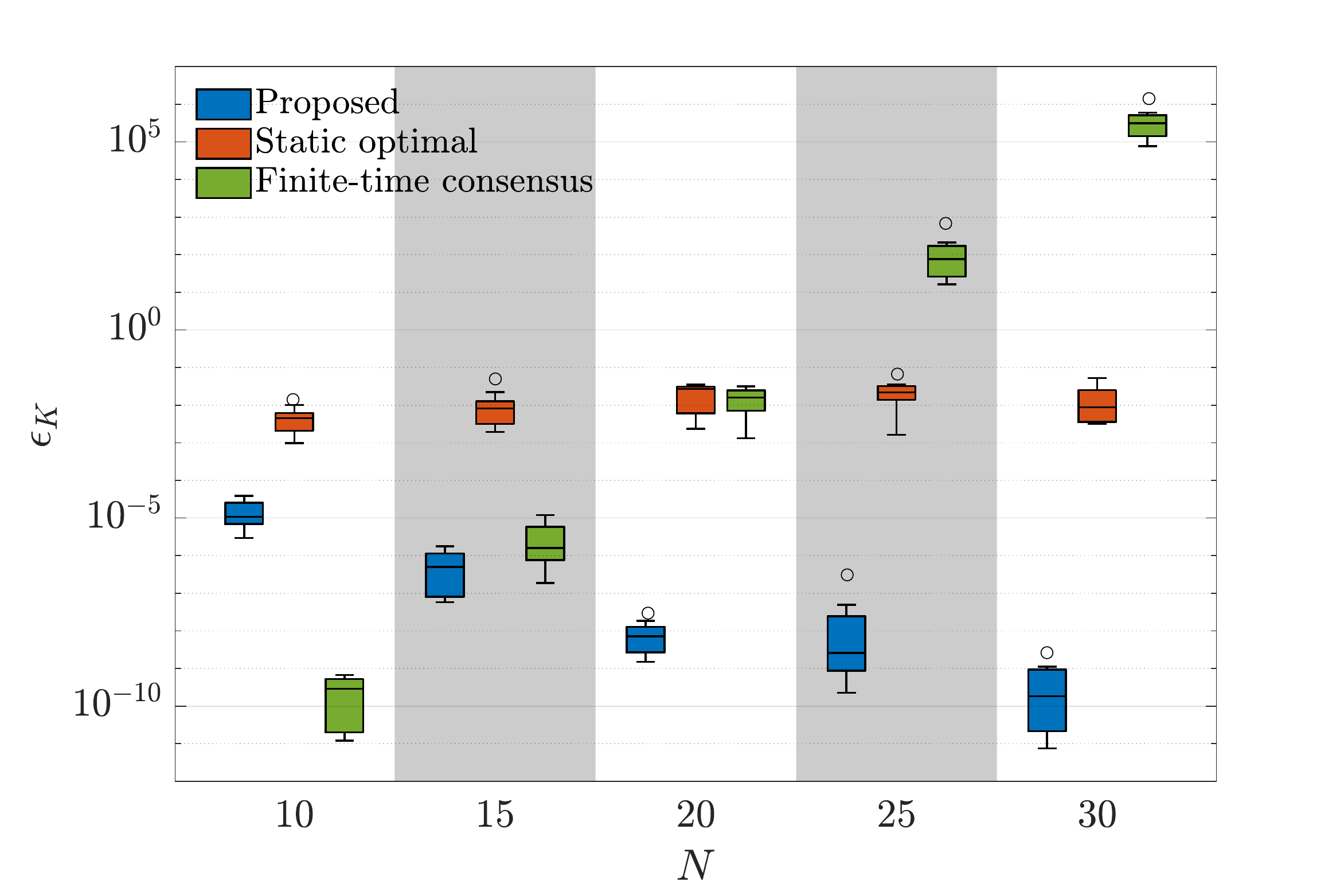}
    \caption{Watts-Strogatz network}
  \end{subfigure}
  \caption{Average consensus errors $\epsilon_K$ for randomly generated small networks.}
  \label{fig:randomGraphs}
\end{figure}

\subsection{Random synthetic networks}

We consider the following three random synthetic network models. The first is the Erd\H{o}s-R\'enyi (ER) network, where we set the probability for edge creation to $0.1$. The second is the Barab\'asi-Albert (BA) model~\cite{Barabasi1999}, where we set the number of edges to attach from a new node to existing nodes to $3$. The third is the Watts-Strogatz (WS) model~\cite{Watts1998}, where each node is joined to its $4$ nearest neighbors in a ring topology and the probability of rewiring each edge is $0.15$.  As in the case of the deterministic networks, we assume that the initial states of the nodes independently follow a uniform distribution on the interval~$[-1, 1]$.

\subsubsection{Small network consensus errors for a finite-time window of length $K$}

We first performed numerical experiments on small-scale networks. We generated 10 networks for each of the three network models with network sizes of $N \in \{10, 15, \dotsc, 30\}$. For each of the generated networks, we optimized the edge weights at the times $k=0, \dotsc, 9$ by minimizing $\epsilon_T$ with $T=10$ and periodically continue the obtained sequence of weighted networks.
We then evaluated the empirical average of the consensus error $\epsilon_K$. We present the results in Fig.~\ref{fig:randomGraphs}. One can observe the same trend as that in Table~\ref{table:}. Specifically, although the finite-time consensus algorithm performs quite well for $N=10$, its consensus error grows exponentially as network size increases, regardless of the network model. In contrast, the accuracy of the consensus achieved by the proposed method is robust to changes in network size and is greater than that of the consensus achieved by the static optimal strategy.

\begin{figure}[tb]
  \centering
  \includegraphics[width=.9\linewidth]{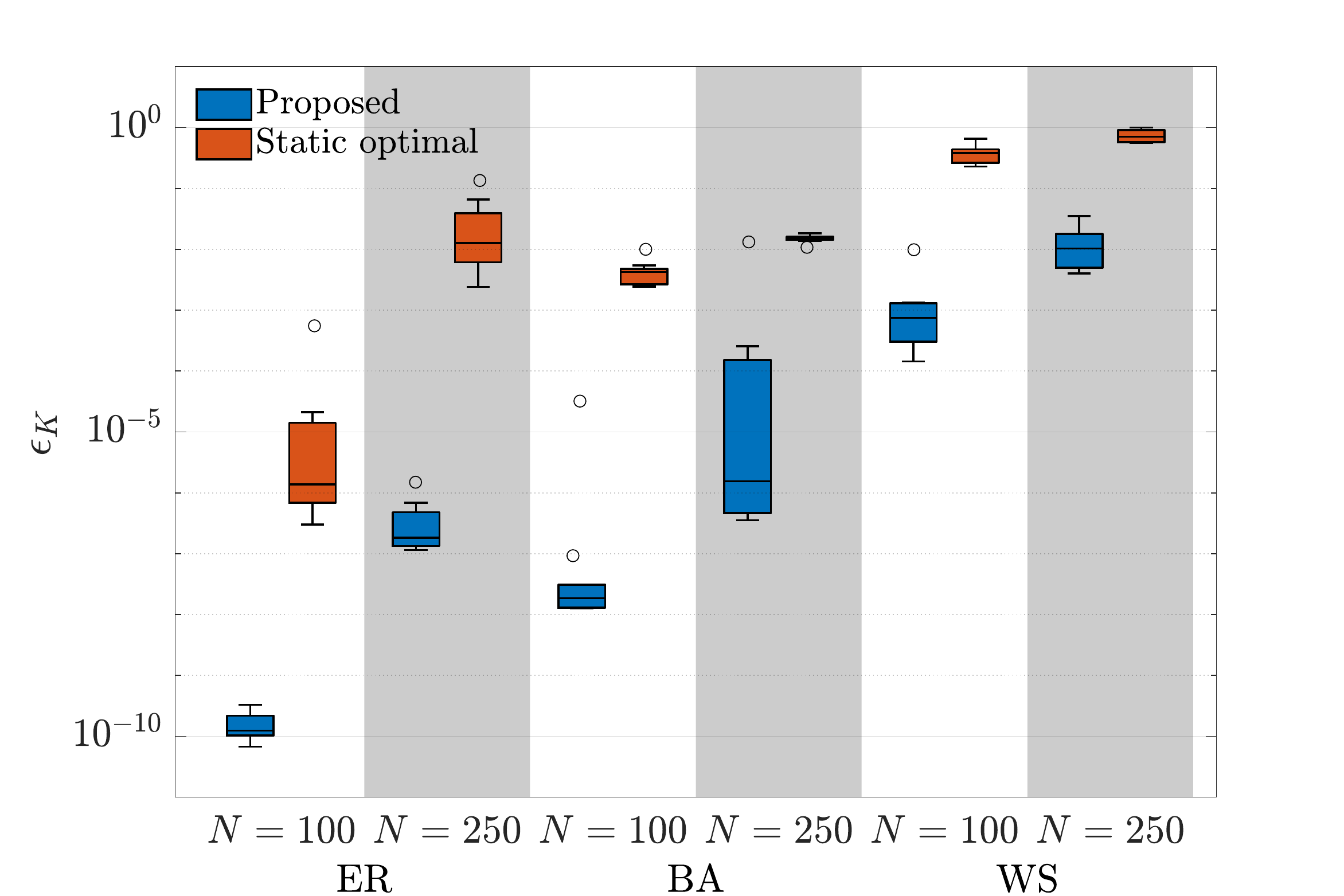}
  \caption{Average consensus errors $\epsilon_K$ for randomly generated networks.}
  \label{fig:largeRandomGraphs}
\end{figure}

\begin{figure}[tb]
  \centering
  \includegraphics[width=.9\linewidth]{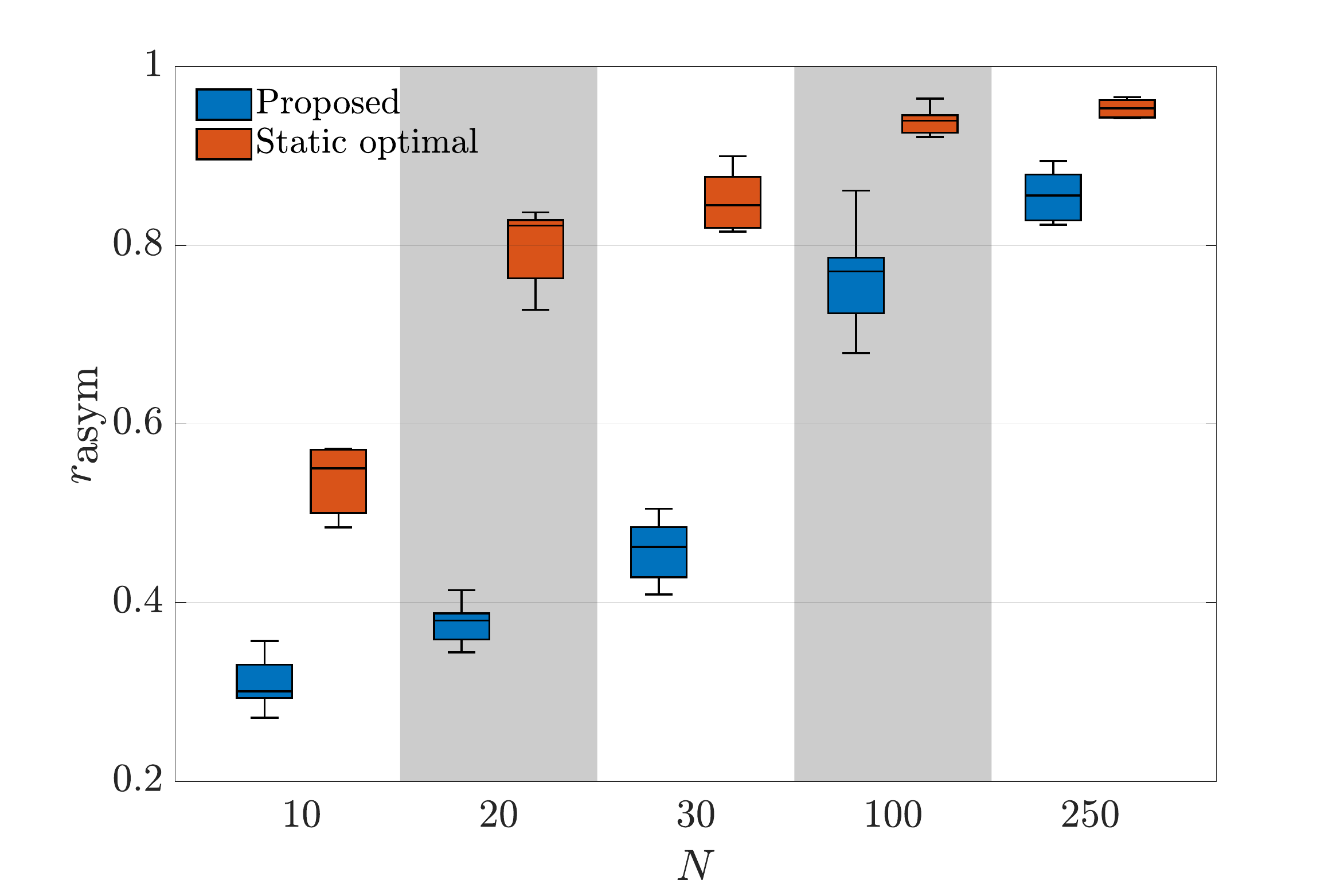}
  \caption{Convergence factors for WS networks.}
  \label{fig:convFactors}
\end{figure}

\subsubsection{Large network consensus errors}

We also conducted numerical experiments on larger networks. We used network sizes of $N \in \{100, 250\}$ and performed experiments similar to those performed for small networks. The finite-time consensus algorithm was excluded because this algorithm generates extremely large consensus errors  due to the numerical instability discussed above. We also reduced the probability for edge creation in the ER network to $0.025$. We made this change because otherwise the linear matrix inequality arising from the static optimal strategy was infeasible. We show the empirical averages of the consensus errors by the proposed method and the static optimal strategy in Fig.~\ref{fig:largeRandomGraphs}. As in Fig.~\ref{fig:randomGraphs} for small networks, one can see that the proposed method yields consensuses with smaller errors compared to the static optimal strategy, regardless of network sizes or models. We can also confirm that the proposed method achieves smaller asymptotic convergence factors. In Fig.~\ref{fig:convFactors}, we present the asymptotic convergence factors for WS networks of various sizes.

\section{Conclusion} \label{sec:conc}

In this paper, we presented a data-driven approach to accelerating the linear average consensus algorithm for undirected temporal networks. The proposed approach first unfolds the consensus algorithm to obtain an equivalent feedforward signal-flow graph, which is regarded as a neural network.  Standard deep learning techniques are then applied to train the obtained neural network, which is a temporal network with optimized edge weights. Numerical experiments confirmed that the proposed method can significantly accelerate the average consensus algorithm for both finite and infinite-time windows.

% \bibliography{IEEEabrv,oguraLibrary2,kishidaLibrary}

\begin{thebibliography}{10}
  \providecommand{\url}[1]{#1}
  \csname url@samestyle\endcsname
  \providecommand{\newblock}{\relax}
  \providecommand{\bibinfo}[2]{#2}
  \providecommand{\BIBentrySTDinterwordspacing}{\spaceskip=0pt\relax}
  \providecommand{\BIBentryALTinterwordstretchfactor}{4}
  \providecommand{\BIBentryALTinterwordspacing}{\spaceskip=\fontdimen2\font plus
  \BIBentryALTinterwordstretchfactor\fontdimen3\font minus
    \fontdimen4\font\relax}
  \providecommand{\BIBforeignlanguage}[2]{{%
  \expandafter\ifx\csname l@#1\endcsname\relax
  \typeout{** WARNING: IEEEtran.bst: No hyphenation pattern has been}%
  \typeout{** loaded for the language `#1'. Using the pattern for}%
  \typeout{** the default language instead.}%
  \else
  \language=\csname l@#1\endcsname
  \fi
  #2}}
  \providecommand{\BIBdecl}{\relax}
  \BIBdecl
  
  \bibitem{Olfati-Saber2007}
  R.~Olfati-Saber, J.~A. Fax, and R.~M. Murray, ``{Consensus and cooperation in
    networked multi-agent systems},'' \emph{Proceedings of the IEEE}, vol.~95,
    no.~1, pp. 215--233, 2007.
  
  \bibitem{Cyb89}
  G.~Cybenko, ``Dynamic load balancing for distributed memory multiprocessors,''
    \emph{Journal of Parallel and Distributed Computing}, vol.~7, no.~2, pp. 279
    -- 301, 1989.
  
  \bibitem{XiaBL05}
  L.~Xiao, S.~Boyd, and S.~Lall, ``A scheme for robust distributed sensor fusion
    based on average consensus,'' in \emph{Fourth International Symposium on
    Information Processing in Sensor Networks}, 2005, pp. 63--70.
  
  \bibitem{Ren2005}
  W.~Ren and R.~W. Beard, ``{Consensus seeking in multiagent systems under
    dynamically changing interaction topologies},'' \emph{IEEE Transactions on
    Automatic Control}, vol.~50, no.~5, pp. 655--661, 2005.
  
  \bibitem{Xiao2004}
  L.~Xiao and S.~Boyd, ``{Fast linear iterations for distributed averaging},''
    \emph{Systems \& Control Letters}, vol.~53, no.~1, pp. 65--78, 2004.
  
  \bibitem{Cortes2005}
  J.~Cort\'es and F.~Bullo, ``{Coordination and geometric optimization via
    distributed dynamical systems},'' \emph{SIAM Journal on Control and
    Optimization}, vol.~44, no.~5, pp. 1543--1574, 2005.
  
  \bibitem{Senel2017}
  K.~Senel and M.~Akar, ``{A distributed coverage adjustment algorithm for
    femtocell networks},'' \emph{IEEE Transactions on Vehicular Technology},
    vol.~66, no.~2, pp. 1739--1747, 2017.
  
  \bibitem{Dorfler2012}
  F.~D\"orfler and F.~Bullo, ``{Synchronization and transient stability in power
    networks and nonuniform Kuramoto oscillators},'' \emph{SIAM Journal on
    Control and Optimization}, vol.~50, no.~3, pp. 1616--1642, 2012.
  
  \bibitem{CheS12}
  J.~Chen and A.~H. Sayed, ``Diffusion adaptation strategies for distributed
    optimization and learning over networks,'' \emph{IEEE Transactions on Signal
    Processing}, vol.~60, no.~8, pp. 4289--4305, 2012.
  
  \bibitem{TsiLR12}
  K.~I. Tsianos, S.~Lawlor, and M.~G. Rabbat, ``Consensus-based distributed
    optimization: Practical issues and applications in large-scale machine
    learning,'' in \emph{Annual Allerton Conference on Communication, Control,
    and Computing}, 2012, pp. 1543--1550.
  
  \bibitem{Olfati-Saber2004}
  R.~Olfati-Saber and R.~M. Murray, ``{Consensus problems in networks of agents
    with switching topology and time-delays},'' \emph{IEEE Transactions on
    Automatic Control}, vol.~49, no.~9, pp. 1520--1533, 2004.
  
  \bibitem{Kempton2018}
  L.~Kempton, G.~Herrmann, and M.~{Di Bernardo}, ``{Self-organization of weighted
    networks for optimal synchronizability},'' \emph{IEEE Transactions on Control
    of Network Systems}, vol.~5, no.~4, pp. 1541--1550, 2018.
  
  \bibitem{Zelazo2013}
  D.~Zelazo, S.~Schuler, and F.~Allg{\"{o}}wer, ``{Performance and design of
    cycles in consensus networks},'' \emph{Systems and Control Letters}, vol.~62,
    no.~1, pp. 85--96, 2013.
  
  \bibitem{Hao2012}
  H.~Hao and P.~Barooah, ``{Improving convergence rate of distributed consensus
    through asymmetric weights},'' in \emph{American Control Conference}, 2012,
    pp. 787--792.
  
  \bibitem{BoyGP06}
  S.~Boyd, A.~Ghosh, B.~Prabhakar, and D.~Shah, ``Randomized gossip algorithms,''
    \emph{IEEE Transactions on Information Theory}, vol.~52, no.~6, pp.
    2508--2530, 2006.
  
  \bibitem{Sundaram2007}
  S.~Sundaram and C.~N. Hadjicostis, ``{Finite-time distributed consensus in
    graphs with time-invariant topologies},'' in \emph{American Control
    Conference}, 2007, pp. 711--716.
  
  \bibitem{SanKM14}
  A.~Sandryhaila, S.~Kar, and J.~M.~F. Moura, ``Finite-time distributed consensus
    through graph filters,'' in \emph{IEEE International Conference on Acoustics,
    Speech, and Signal Processing}, 2014, pp. 1080--1084.
  
  \bibitem{Hendrickx2015}
  J.~M. Hendrickx, G.~Shi, and K.~H. Johansson, ``{Finite-time consensus using
    stochastic matrices with positive diagonals},'' \emph{IEEE Transactions on
    Automatic Control}, vol.~60, no.~4, pp. 1070--1073, 2015.
  
  \bibitem{Safavi2015}
  S.~Safavi and U.~A. Khan, ``{Revisiting finite-time distributed algorithms via
    successive nulling of eigenvalues},'' \emph{IEEE Signal Processing Letters},
    vol.~22, no.~1, pp. 54--57, 2015.
  
  \bibitem{Shang2016}
  Y.~Shang, ``{Finite-time weighted average consensus and generalized consensus
    over a subset},'' \emph{IEEE Access}, vol.~4, no.~8, pp. 2615--2620, 2016.
  
  \bibitem{Shuman2013}
  D.~I. Shuman, S.~K. Narang, P.~Frossard, A.~Ortega, and P.~Vandergheynst,
    ``{The emerging field of signal processing on graphs: Extending
    high-dimensional data analysis to networks and other irregular domains},''
    \emph{IEEE Signal Processing Magazine}, vol.~30, no.~3, pp. 83--98, 2013.
  
  \bibitem{Apers2017}
  S.~Apers and A.~Sarlette, ``{Accelerating consensus by spectral clustering and
    polynomial filters},'' \emph{IEEE Transactions on Control of Network
    Systems}, vol.~4, no.~3, pp. 544--554, 2017.
  
  \bibitem{Falsone2018}
  A.~Falsone, K.~Margellos, S.~Garatti, and M.~Prandini, ``{Finite-time
    distributed averaging over gossip-constrained ring networks},'' \emph{IEEE
    Transactions on Control of Network Systems}, vol.~5, no.~3, pp. 879--887,
    2018.
  
  \bibitem{Ito2019}
  D.~Ito, S.~Takabe, and T.~Wadayama, ``{Trainable ISTA for sparse signal
    recovery},'' \emph{IEEE Transactions on Signal Processing}, vol.~67, no.~12,
    pp. 3113--3125, 2019.
  
  \bibitem{pytorch}
  A.~Paszke, S.~Gross, S.~Chintala, G.~Chanan, E.~Yang, Z.~DeVito, Z.~Lin,
    A.~Desmaison, L.~Antiga, and A.~Lerer, ``Automatic differentiation in
    pytorch,'' in \emph{31st Conference on Neural Information Processing
    Systems}, 2017.
  
  \bibitem{Boyd1994}
  S.~Boyd, L.~{El Ghaoui}, E.~Feron, and V.~Balakrishnan, \emph{{Linear Matrix
    Inequalities in System and Control Theory}}.\hskip 1em plus 0.5em minus
    0.4em\relax Society for Industrial Mathematics, 1994.
  
  \bibitem{Zachary1977}
  W.~W. Zachary, ``{An information flow model for conflict and fission in small
    groups},'' \emph{Journal of Anthropological Research}, vol.~33, no.~4, pp.
    452--473, 1977.
  
  \bibitem{Barabasi1999}
  A.-L. Barab\'asi and R.~Albert, ``{Emergence of scaling in random networks},''
    \emph{Science}, vol. 286, no. 5439, pp. 509--512, 1999.
  
  \bibitem{Watts1998}
  D.~J. Watts and S.~H. Strogatz, ``{Collective dynamics of 'small-world'
    networks},'' \emph{Nature}, vol. 393, no. 6684, pp. 440--442, 1998.
  
  \end{thebibliography}


\appendix

\section{Proof of Proposition~\ref{LEM:}}

In this appendix, we present the proof of Proposition~\ref{LEM:}. 
We begin by presenting a few lemmas. For a real sequence $a  =\{a(k)\}_{k=0}^\infty$, we define
\begin{equation}\label{eq:def:eta}
  \eta(a) = \limsup_{k\to\infty} \abs{a(k)}^{1/k} \in [0, \infty].
\end{equation}
It should be noted that the following relationship holds: 
\begin{equation}\label{eq:logeta}
  \log \eta(a) = \limsup_{k\to \infty} \frac{\log \abs{a(k)}}{k}.
\end{equation}

\begin{lemma}\label{lem:eta:CL}  
Let $a  =\{a(k)\}_{k=0}^\infty$ and $b = \{b(k)\}_{k=0}^\infty$ be real sequences. Assume that there exist integers~$L_1$, $L_2$ and positive constants $C_1$, $C_2$ such that
\begin{equation}\label{eq:C1C2L1L2}
  C_1 \abs{b(k+L_1)} \leq \abs{a(k)} \leq C_2 \abs{b(k+L_2)}
\end{equation}
for all $k\geq \max(0, -L_1, -L_2)$. Then, $\eta(a) = \eta(b)$. 
\end{lemma}

\begin{IEEEproof}
By taking the logarithms in the inequality~\eqref{eq:C1C2L1L2}, we obtain
    \begin{equation*}
      \begin{multlined}
        \frac{\log C_1}{k} + \frac{\log \abs{b(k+L_1)}}{k+L_1}\cdot \frac{k+L_1}{k}
        \leq
        \frac{\log \abs{a(k)}}{k}
        \\
        \leq
        \frac{\log C_2}{k} + \frac{\log \abs{b(k+L_2)}}{k+L_2}\cdot \frac{k+L_2}{k}
      \end{multlined}
    \end{equation*}
    for all $k \geq \max(0, -L_1, -L_2)$.
As desired,  taking the limit superiors with respect to $k$ in this inequality and using \eqref{eq:logeta} show $\log \eta(a) = \log \eta(b)$.
\end{IEEEproof}

\begin{lemma}\label{lem:eta:T} 
  Let $a  =\{a(k)\}_{k=0}^\infty$ be a real sequence. Let $T$ be a positive integer. For a nonnegative integer~$k$, we define 
    \begin{equation*}
      \lfloor k\rfloor_T = \lfloor k/T \rfloor T,  
    \end{equation*}
    where $\lfloor \cdot \rfloor$ denotes the floor function. If $b(k) = a(\lfloor k \rfloor_T)$ for all $k\geq 0$, then $\eta(b) = \eta(\{a(kT) \}_{k\geq 0})^{1/T}$.
\end{lemma}
\begin{IEEEproof}
  The sequence $\{ \sup_{\ell\geq k} {\ell}^{-1}{\log \abs{b(\ell)}} \}_{k=1}^\infty$ is convergent in the extended real space $[-\infty, \infty]$.
  Therefore, its subsequence $\{ \sup_{\ell\geq kT} {\ell}^{-1}{\log \abs{b(\ell)}} \}_{k=1}^\infty$ is convergent and converges to the limit of the original sequence. Hence, from \eqref{eq:logeta}, we can obtain
  \begin{equation}\label{eq:logetab1}
    \begin{aligned}
      \log \eta(b)
       & =
       \adjustlimits \lim_{k\to\infty} \sup_{\ell\geq k} \frac{\log \abs{b(\ell)}}{\ell}
      \\
       & =
       \adjustlimits \lim_{k\to\infty} \sup_{\ell\geq kT} \frac{\log \abs{b(\ell)}}{\ell}
      \\
       & =
       \adjustlimits \lim_{k\to\infty} \sup_{k'\geq k} \max_{k'T \leq \ell < (k'+1)T} \frac{\log \abs{b(\ell)}}{\ell}.
    \end{aligned}
  \end{equation}
  Because $b(\ell) = a(k' T)$ for any integer~$\ell$ satisfying $k'T \leq \ell < (k'+1)T$, we have
  \begin{equation}\label{eq:logetab2}
    \begin{aligned}
      \max_{k'T \leq \ell < (k'+1)T } \frac{\log \abs{b(\ell)}}{\ell}
       & =
      \max_{k'T \leq \ell < (k'+1)T } \frac{\log \abs{a(k' T)}}{\ell}
      \\
       & =
      \frac{\log \abs{a(k' T)}}{k' T}.
    \end{aligned}
  \end{equation}
As desired, equations~\eqref{eq:logetab1} and~\eqref{eq:logetab2} yield 
  \begin{equation*}
    \log \eta(b) = \adjustlimits \lim_{k\to\infty} \sup_{k'\geq k} \frac{\log \abs{a(k' T)}}{k' T}
    = \frac{1}{T}\log \eta(\{a(kT) \}_{k\geq 0}).
  \end{equation*}
\end{IEEEproof}

\begin{lemma}\label{lem:eta:M} 
  Let $M$ be an ${n\times n}$ real matrix with generalized eigenvectors $v_1$, \dots, $v_n \in \mathbb{R}^n$ corresponding to the eigenvalues~$\lambda_1$, \dots, $\lambda_n$, counted according to algebraic multiplicity. Let $x_0 \in \mathbb{R}^n$ be arbitrary and assume that there exist a set of integers $\mathcal I \subset\{1, \dotsc, n\}$ and nonzero numbers $c_i$ ($i\in\mathcal I$) such that $x_0 = \sum_{i\in \mathcal I} c_i v_i$. Then, the sequence $a = \{ \norm{ M^k x_0}\}_{k=0}^\infty$ satisfies $\eta(a) = \max_{i\in\mathcal I}\abs{\lambda_i}$. 
\end{lemma}
\begin{IEEEproof}
 Let $p$ denote the number of distinct eigenvalues of $M$.
    Without loss of generality, we can assume that the matrix~$M$ is in the Jordan canonical form with the Jordan blocks 
    \begin{equation*}
      J_k = c_k I + N_k\in \mathbb{R}^{n_k \times n_k},\ k=1, \dotsc, p, 
    \end{equation*}
    with a real number~$c_k$, a nilpotent matrix~$N_k$, and a positive integer $n_k$ satisfying $\sum_{k=1}^pn_k = n$. We can also assume the existence of a positive integer $q \leq p$ such that each of the generalized eigenvectors~$v_{i}$ ($i\in\mathcal I$) corresponds to one of the first $q$ Jordan blocks $J_1$, \dots, $J_q$. This assumption implies the following identity:
 \begin{equation}\label{eq:clambda}
      \max_{1\leq k\leq q} \abs{c_k} = \max_{i\in\mathcal I}\abs{\lambda_i},
    \end{equation}
as well as the existence of nonzero vectors $\xi_k \in \mathbb{R}^{n_k}$ ($k=1, \dotsc, q$) such that
    \begin{equation*}
      x_0 = \begin{bmatrix}
        \xi_1 \\ \vdots \\ \xi_q \\ 0_{n-\sum_{k=1}^q n_k}
      \end{bmatrix}, 
    \end{equation*}
    where $0_{n-\sum_{k=1}^q n_k}$ denotes the zero vector of length $n-\sum_{k=1}^q n_k$.
    Therefore, we obtain 
    \begin{equation*}
      M^k x_0 = \begin{bmatrix}
        J_1^k \xi_1 \\ \vdots \\ J_q^k \xi_q \\ 0_{n-\sum_{k=1}^q n_k}
      \end{bmatrix}
      =
      \begin{bmatrix}
        c_1^k w_{1}(k) \\ \vdots \\ c_q^k w_{q}(k) \\ 0_{n-\sum_{k=1}^q n_k}
      \end{bmatrix}, 
    \end{equation*}
    where $w_{1}$, \dots, $w_{q}$ are nonzero vectors growing polynomially in $k$.
    Hence, the definition of $\eta$ yields $\eta(a) =  \max_{1\leq k\leq q} \abs{c_k}$. 
    This equation and \eqref{eq:clambda} complete the proof.
\end{IEEEproof}

We are now ready to prove Proposition~\ref{LEM:}.
    Using the notation~\eqref{eq:def:eta}, we obtain
    \begin{equation}\label{eq:rasym:lemma}
      r_{\asym} = \sup_{x_0 \neq c\onev} \eta(\{ \norm{e(k)} \}_{k=0}^\infty)
    \end{equation}
    because $\norm{e(0)}^{1/k}$ converges to one as $k\to\infty$. Equation~\eqref{eq:periodicConsensnsProt} shows that the error vector~$e$ satisfies
    \begin{align*}\begin{aligned}
        e(sT+\tau+1) & =\left(\prod_{t=0}^{\tau}(I-L^\star(\tau-t)) \right)e(sT)
      \end{aligned}
    \end{align*}
    for all $0\leq \tau \leq T-1$ and $s\geq 0$. Therefore, if we define
    \begin{equation*}
      C = \max_{0\leq \tau \leq T-1} \prod_{t=0}^{\tau} \norm{I-L^\star(t)},
    \end{equation*}
    then we can show that $C^{-1} \norm{e(\lfloor k+T\rfloor_T)} \leq \norm{e(k)} \leq C \norm{e(\lfloor k\rfloor_T)}$ for all $k\geq 0$. By Lemmas~\ref{lem:eta:CL} and~\ref{lem:eta:T}, we obtain
    \begin{equation}\label{eq:normek}
      \begin{aligned}
        \eta(\{ \norm{e(k)} \}_{k=0}^\infty)
         & =
        \eta(\{ \norm{e(\lfloor k \rfloor_T)} \}_{k=0}^\infty)
        \\
         & =
        \eta(\{ \norm{e(kT)} \}_{k=0}^\infty)^{1/T}
        \\
         & =
        \eta(\{ \norm{\bar e(k)} \}_{k=0}^\infty)^{1/T},
      \end{aligned}
    \end{equation}
    where $\bar e$ is defined as $\bar e(k) = e(kT)$ for all $k\geq 0$.

  It should be noted that the sequence~$\{\bar e(k)\}_{k=0}^\infty$ satisfies
  $\bar e(k+1) = M\bar e(k)$
  for all $k\geq 0$.
  Additionally, based on the assumption in the lemma, the eigenvalue~$1$ of $M$ corresponding to the eigenvector~$\onev$ is simple.
  Because $\bar e(0) = e(0)= x_0-c\onev$ belongs to the space spanned by the generalized eigenvectors of $M$ corresponding to the other eigenvalues, we have $\eta(\{ \norm{\bar e(k)} \}_{k=0}^\infty) \leq \rho$ by Lemma~\ref{lem:eta:M}. The equality is attained when $x_0$ equals one of the generalized eigenvectors of $M$ corresponding to the eigenvalue having the modulus $\rho$. Therefore, equations~\eqref{eq:rasym:lemma} and~\eqref{eq:normek} complete the proof. 

\bibliographystyle{IEEEtran}
\begin{IEEEbiography}[{\includegraphics[width=1in,height=1.25in,clip,keepaspectratio]{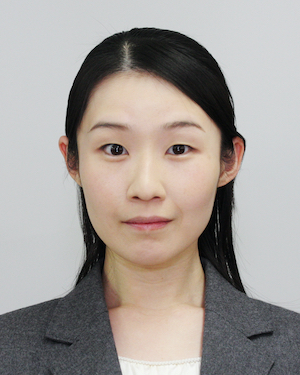}}]{Masako Kishida}(M'10--SM'18) received the Ph.D. degree in Mechanical Engineering from University of Illinois at Urbana-Champaign, IL, in 2010. From 2010 to 2016, she held positions in the U.S.A., Japan, New Zealand and Germany.
  Since 2016, she has been an Associate Professor at National Institute of Informatics, Tokyo Japan.
  Her research interests include networked control systems, optimizations and uncertainty analysis.  Prof. Kishida was a recipient of the Humboldt Research Fellowship in 2015, the Telecom System Technology Award in 2019 and the Young Scientists' Prize for the Commendation of Science and Technology by MEXT in 2020. 
\end{IEEEbiography}

\begin{IEEEbiography}[{\includegraphics[width=1in,height=1.25in,clip,keepaspectratio]{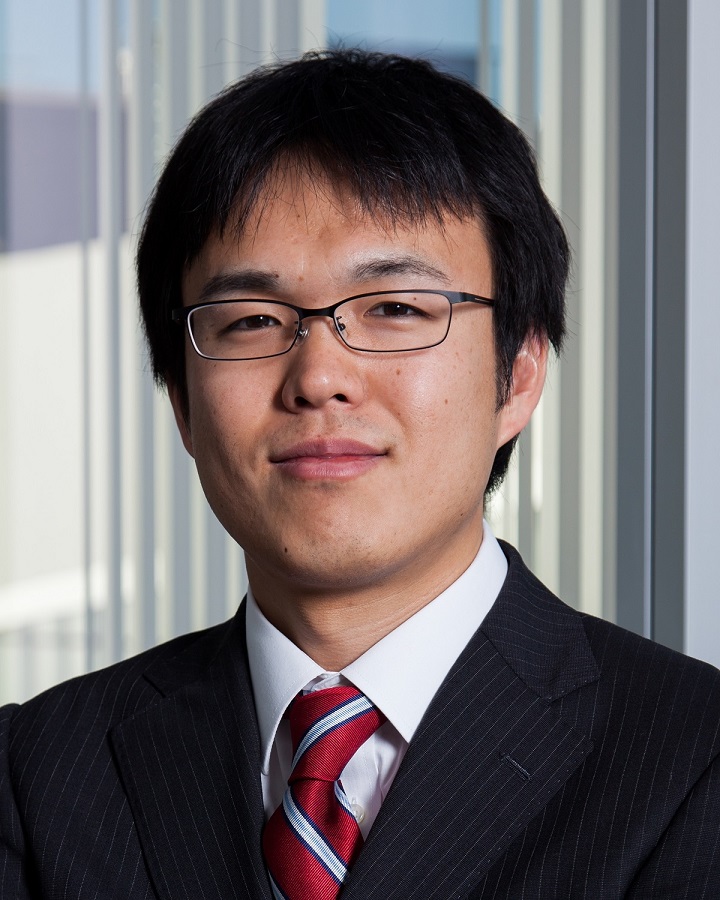}}]{Masaki Ogura}
  (M'14) is an Associate Professor in the Graduate School of Information Science and Technology at Osaka University, Japan. Prior to joining Osaka University, he was a Postdoctoral Researcher at the University of Pennsylvania, USA and an Assistant Professor at the Nara Institute of Science and Technology, Japan. His research interests include network science, dynamical systems, and stochastic processes with applications in networked epidemiology, design engineering, and biological physics. He was a runner-up of the 2019 Best Paper Award by the IEEE Transactions on Network Science and Engineering and a recipient of the 2012 SICE Best Paper Award. He is an Associate Editor of the Journal of the Franklin Institute.
\end{IEEEbiography}

\begin{IEEEbiography}[{\includegraphics[width=1in,height=1.25in,clip,keepaspectratio]{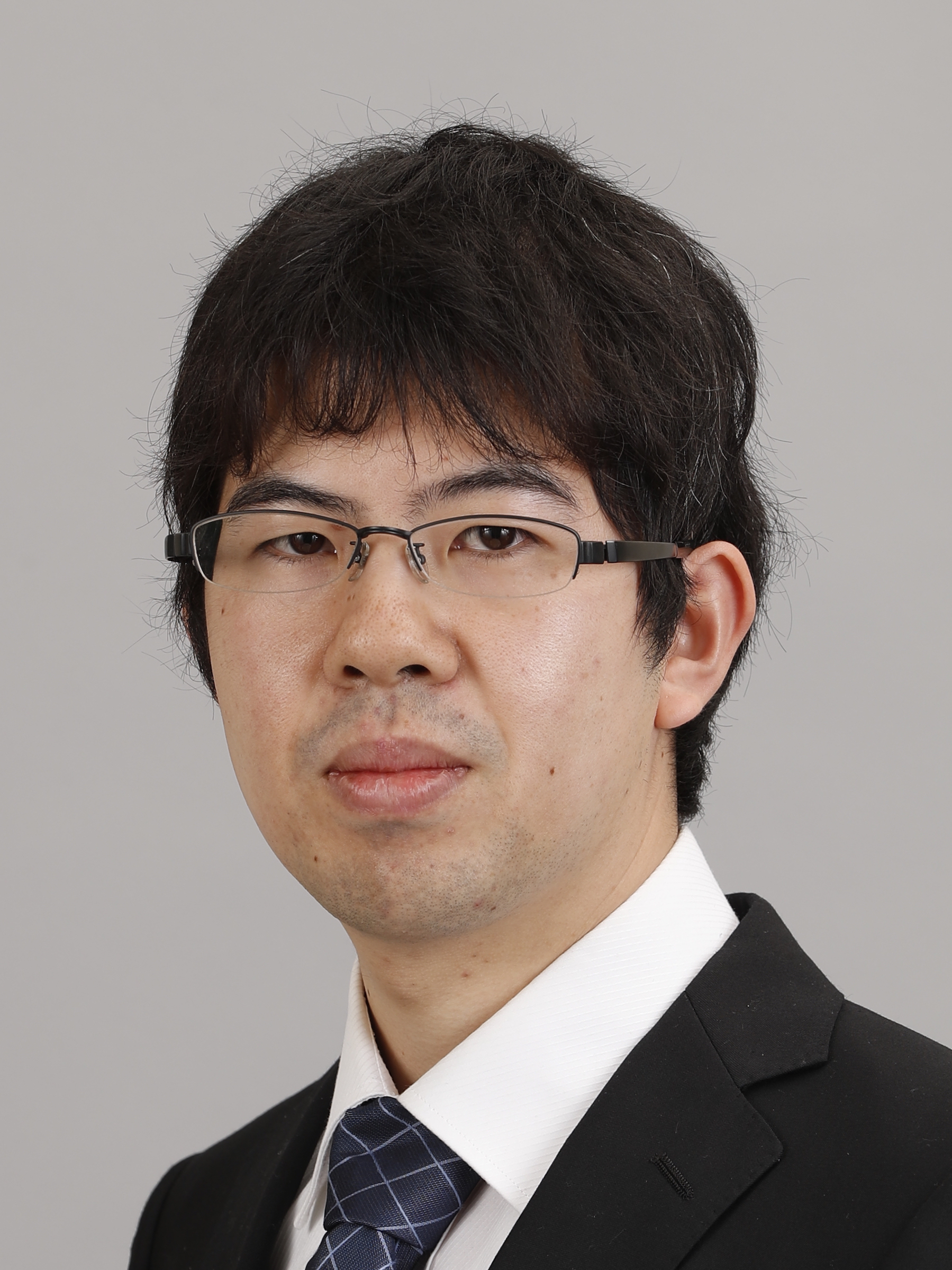}}]{Yuichi Yoshida} received the B.S. degree in Engineering and the M.S. and Ph.D. degrees in Informatics, all from Kyoto University, in 2007, 2009, and 2012, respectively. In 2012, he joined National Institute of Informatics, Tokyo, Japan as an assistant professor, and since 2015, he has been an associate professor there. Prof. Yoshida received Kyoto University President Award in 2012, JSPS Ikushi Prize in 2012, the Young Scientists’ Prize for the Commendation of Science and Technology by MEXT in 2017, IPSJ Microsoft Research Award on Information in 2018, and AISTATS Best Paper Award in 2018.
\end{IEEEbiography}

\begin{IEEEbiography}[{\includegraphics[width=1in,height=1.25in,clip,keepaspectratio]{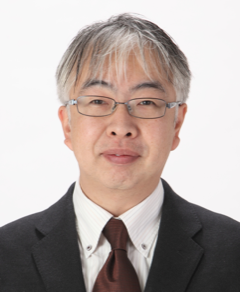}}]{Tadashi Wadayama}
(M'96) was born in Kyoto, Japan, on May 9,1968. 
He received the B.E., the M.E., and the D.E. degrees from Kyoto Institute of Technology in 1991, 1993 and 1997, respectively. 
On 1995, he started to work with Faculty of Computer Science and System Engineering, Okayama Prefectural University as a research associate. 
From April 1999 to March 2000, he stayed in Institute of Experimental Mathematics, Essen University (Germany) as a visiting researcher. 
On 2004, he moved to Nagoya Institute of Technology as an associate professor. Since 2010, he has been a full professor of Nagoya Institute of Technology. 
His research interests are in coding theory, information theory, and coding and signal processing for digital communication/storage systems. 
He is a member of IEICE. 
\end{IEEEbiography}
\EOD

\end{document}